\documentclass[11pt,a4paper]{article}

\usepackage{inputenc}
\usepackage{amsmath}
\usepackage{bm}
\usepackage{bbold}
\usepackage{amsthm}
\usepackage{enumerate}

\usepackage{algorithm}

\usepackage[top=1.0in, left=1.0in, right=1.0in, bottom=1.0in]{geometry}

\usepackage{tikz}\tikzset{x=1cm,y=1cm,z=1cm}

\usepackage{pgfplots}\pgfplotsset{compat=1.16}

\usepackage[hyphens]{url}

\usepackage{hyperref}
\usepackage{breakurl}

\title{Tropical solution of discrete best approximation problems\thanks{Mathematics, 2025, 13(22), 3660;  https://doi.org/10.3390/math13223660}}

\author{N. Krivulin\thanks{Faculty of Mathematics and Mechanics, St.~Petersburg State University, 28 Universitetsky Ave., St.~Petersburg, 198504, Russia; 
nkk@math.spbu.ru.}
}

\date{}

\newtheorem{theorem}{Theorem}
\newtheorem{lemma}[theorem]{Lemma}

\theoremstyle{definition}

\begin{document}

\maketitle

\begin{abstract}
We consider discrete best approximation problems in the setting of tropical algebra, which is concerned with the theory and application of algebraic systems with idempotent operations. Given a set of input--output pairs of an unknown function defined on a tropical semifield, the problem is to determine an approximating rational function formed by two Puiseux polynomials as numerator and denominator. With specified numbers of monomials in both polynomials, the approximation aims at evaluating the exponent and coefficient for each monomial in the polynomials to fit the rational function to the data in the sense of a tropical distance function. To solve the problem, we transform it into an approximation of a vector equation with unknown vectors on both sides, where one side corresponds to the numerator polynomial and the other side to the denominator. Each side involves a matrix with entries dependent on the unknown exponents, multiplied by the vector of unknown coefficients of monomials. We propose an algorithm that constructs a series of approximate solutions by alternately fixing one side of the equation to an already-found result and leaving the other side intact. Each equation obtained is approximated with respect to the vector of coefficients, which yields this vector and approximation error, both parameterized by  exponents. The exponents are found by minimizing the error with an optimization procedure based on an agglomerative clustering technique. To illustrate, we present results for an approximation problem in terms of max-plus algebra (a real semifield with addition defined as maximum and multiplication as arithmetic addition), which corresponds to an ordinary problem of piecewise linear approximation of real functions. As our numerical experience shows, the proposed algorithm converges in a finite number of steps and provides a reasonably accurate solution to the problems considered.
\\

\textbf{Key-Words:} tropical semifield, tropical Puiseux polynomial, best approximate solution, discrete best approximation, Chebyshev approximation.
\\

\textbf{MSC (2020):} 15A80, 90C24, 41A50, 41A65, 65D15
\end{abstract}

\section{Introduction}

Discrete best approximation problems, where sample data are fitted to a parameterized function by adjusting parameters, constitute an essential area of study in approximation theory and its applications \cite{Mhaskar2000Fundamentals,Steffens2006History}. The best approximation problems are formulated to find those values of the parameters which provide a minimum discrepancy between the sample data and approximating function in the sense of some metric defined on a vector space to represent the approximation error. As notable examples of approximation problems, which date back to Laplace’s classic work \cite{Laplace1832Mecanique} (Book 3, Chapter V, \S39), one can consider problems where the approximating functions are piece-wise linear and the error is measured by the Chebyshev metric. These problems are still of sufficient interest and find application in the solution of overdetermined systems of linear equations, regression analysis when the distribution of errors has bounded support, and  other research and applied contexts.

Since the 1960s, a variety of solutions have been developed to handle the discrete best approximation problems with both Chebyshev and other metrics, as shown in the overviews in \cite{Conn1988Computational,Szusz2010Linear}. Methods and techniques used to solve the problems include linear programming \cite{Stiefel1960Note,Osborne1967Best,Watson1970Algorithm,Sposito1976Minimizing}, dynamic programming \cite{Gluss1962Further,Camponogara2015Models}, and other optimization approaches~\cite{Stone1961Approximation,Cameron1966Piecewise,Tomek1974Two,Imai1986Optimal,Szusz2010Linear,Dellaccio2025Truncated}. The proposed solutions normally offer optimal or near-optimal results that can be obtained with moderate polynomial computational complexity.

{Another approach to formulate and solve best approximation problems is based on the models and methods of tropical algebra \cite{Kolokoltsov1997Idempotent,Golan2003Semirings,Heidergott2006Maxplus,Itenberg2007Tropical,Gondran2008Graphs,Butkovic2010Maxlinear,Maclagan2015Introduction,Kenoufi2025Idempotent}. This algebra is concerned with the theory of algebraic systems with idempotent operations and finds applications in a range of areas from algebraic geometry to operations research. An example of such systems is max-plus algebra, which is a tropical (idempotent) semifield thought of as a set of reals, where addition is defined as maximum and multiplication as arithmetic addition. Another example is max-algebra, a tropical semifield of nonnegative reals, where addition is defined as maximum and multiplication as usual.}

Reformulation of a range of problems that are nonlinear in conventional mathematics into the setting of tropical algebra transforms them into linear problems in the tropical sense, which facilitates the formal analysis and simplifies the derivation of a solution. Using the representation in terms of tropical algebra appears to be fruitful in problems that use mathematical models based on max or min operations, including minimax optimization problems. Methods of optimization in the context of tropical algebra (tropical optimization) are applied to solve many real-world problems in project scheduling, location analysis, decision-making and other fields. In many cases, the methods of tropical optimization can provide analytical solutions in compact closed form for problems that are known to have only numerical algorithmic solutions available.

Existing solutions to approximation problems in the tropical algebra setting mainly concentrate on the best approximation of vectors in tropical linear spaces \cite{Butkovic2010Maxlinear,Akian2011Best,Saad2021Zerosum}. Applications of tropical algebra to solve discrete best approximation problems of functions, which can be found in recent works on neural networks and machine learning (see, e.g., \cite{Zhang2018Tropical,Maragos2021Tropical,Dunbar2024Alternating,Ioannis2025Revisiting}), present another line of investigation that needs further exploration. As an attempt to address this need, a general problem of discrete best approximation with respect to a generalized metric is introduced and examined in the framework of tropical algebra in \cite{Krivulin2023Algebraic,Krivulin2024Solution}. The problem is formulated to approximate functions defined on tropical semifields with generalized tropical Puiseux polynomials and rational functions used as approximants.

A Puiseux polynomial is defined in the same way as classical polynomials in which the exponents can be set to rational numbers. Tropical Puiseux polynomials that appear in the approximation, have the form where addition and multiplication are defined in terms of tropical algebra \cite{Markwig2010Field,Grigoriev2018Tropical}. Unlike the ordinary Puiseux polynomials with exponents assumed to be rational, the generalized polynomials may have real exponents. The Puiseux polynomials find applications in a range of fields, including algebraic geometry \cite{Itenberg2007Tropical,Maclagan2015Introduction}, image processing \cite{Li1992Morphological,Wang2023Tropical}, cryptography \cite{Chen2024Tropical,Durcheva2025Closer}, game theory \cite{Esparza2008Approximative,Bereau2024Tropical}, neural networks \cite{Zhiwei2025Achieving,Ioannis2025Revisiting} and others. The rational functions used as approximants take the form of tropical ratios of Puiseux polynomials. In the context of max-plus algebra, the rational functions define ordinary piecewise linear functions, whereas the tropical metric coincides with the standard Chebyshev metric. This offers a great opportunity to handle conventional discrete best Chebyshev approximation problems through the solution of corresponding approximation problems in the tropical algebra setting and makes the development of new methods and techniques of tropical discrete best approximation problems important. 

In this paper, we consider discrete best approximation problems that are formulated and solved in the framework of tropical algebra. Given a set of input--output pairs of an unknown function defined on a tropical semifield, the problem is to determine an approximating rational function formed by two Puiseux polynomials as its numerator and denominator. With specified numbers of monomials in both polynomials, the approximation aims at evaluating the exponent and coefficient for each monomial in the polynomials to fit the rational function to the given data in the sense of a tropical distance function.

To solve the approximation problem, we follow the approach developed in \cite{Krivulin2023Solution,Krivulin2024Solution} to transform the problem into an approximation of a vector equation with unknown vectors on both sides, where one side corresponds to the numerator polynomial and the other side to the denominator. Each side of the equation involves a matrix with entries dependent on the unknown exponents, multiplied by the vector of unknown coefficients of monomials in the polynomial. We develop a computational procedure based on the alternating algorithm proposed in \cite{Krivulin2023Algebraic}, which constructs a series of approximate solutions by alternately fixing one side of the equation to an already found result and leaving the other intact. 

In the same way as in \cite{Krivulin2024Solution}, each obtained equation with one fixed side is first approximated with respect to the vector of coefficients. We apply the solution proposed in \cite{Krivulin2009Onsolution,Krivulin2012Solution}, which results in a vector of coefficients and approximation error, both parameterized by the exponents. Furthermore, the values of exponents are found by minimizing the approximation error with an optimization procedure developed in \cite{Krivulin2024Solution} on the basis of an agglomerative clustering technique.
To illustrate applications, we present results for approximation problems formulated in terms of max-plus algebra, which correspond to ordinary problems of piecewise linear approximation of real functions.

We consider the purpose of this paper as the first communication about the proposed approximation approach, mainly aimed at a general description of the method and confirmation of its workability. A formal convergence analysis that allows estimating the computational complexity of the solution, as well as analytical and experimental comparison with other existing solutions, is intended to be a subject for further research.

The rest of the paper is organized as follows. We start in Section~\ref{S-PDNR} with an overview of key definitions, notations, and preliminary results of tropical algebra. We continue in Section~\ref{S-TPRF} to consider tropical polynomials and rational functions. Section~\ref{TBDAP} introduces and discusses problems of tropical best discrete approximation by polynomials and rational functions. We describe the procedure of polynomial approximation in detail in Section~\ref{S-APF}. Section~\ref{S-ARF} presents the main outcome of the study, which combines previous results to develop a procedure of rational approximation. We give numerical and graphical examples in Section~\ref{S-NEGI} to illustrate the findings. Section~\ref{S-C} offers some concluding remarks.

\section{Preliminary Definitions, Notations and Results}
\label{S-PDNR}

In this section, we outline basic definitions and notations, and describe preliminary results of tropical algebra to present an analytical framework for solving the approximation problems under consideration. For additional details on the theory and methods of tropical mathematics, one may refer to monographs and textbooks \cite{Kolokoltsov1997Idempotent,Golan2003Semirings,Heidergott2006Maxplus,Itenberg2007Tropical,Gondran2008Graphs,Butkovic2010Maxlinear,Maclagan2015Introduction,Kenoufi2025Idempotent}, and references therein.

\subsection{Tropical Semifield}

Consider a set $\mathbb{X}$ that is closed under operations $\oplus$ (addition) and $\otimes$ (multiplication) and has elements $\mathbb{0}$ (zero) and $\mathbb{1}$ (one). We assume that $(\mathbb{X},\oplus,\mathbb{0})$ is a commutative idempotent monoid, $(\mathbb{X}\setminus\{\mathbb{0}\},\otimes,\mathbb{1})$ is an Abelian group, and multiplication $\otimes$ distributes over addition $\oplus$. Under these assumptions, the algebraic system $(\mathbb{X},\oplus,\otimes,\mathbb{0},\mathbb{1})$ is referred to as the tropical (or idempotent) semifield.

In the tropical semifield, addition is idempotent, which means that $x\oplus x=x$ for each $x\in\mathbb{X}$. The idempotent addition induces on $\mathbb{X}$ a partial order such that $x\leq y$ if and only if $x\oplus y=y$. With respect to this order, addition and multiplication are monotone in the sense that if $x\leq y$ for $x,y\in\mathbb{X}$, then $x\oplus z\leq y\oplus z$ and $x\otimes z\leq y\otimes z$ for any $z\in\mathbb{X}$. Addition possesses an extremal property (majority law) which says that $x\leq x\oplus y$ and $y\leq x\oplus y$. Finally, the inequality $x\oplus y\leq z$ is equivalent to the system of inequalities $x\leq z$ and $y\leq z$.

In what follows, we further assume that addition satisfies the property of selectivity in the form $x\oplus y\in\{x,y\}$, and thus the order associated with the addition is a total order. 

Multiplication is invertible, which provides any $x\ne\mathbb{0}$ with its inverse $x^{-1}$ such that $x\otimes x^{-1}=\mathbb{1}$. The integer powers are defined in the usual way as $x^{0}=\mathbb{1}$, $x^{p}=x\otimes x^{p-1}$, $x^{-p}=(x^{-1})^{p}$ and $\mathbb{0}^{p}=\mathbb{0}$ for any $x\ne\mathbb{0}$ and integer $p>0$. It is assumed that the integer powers can be extended to rational and then to real powers to make exponentiation with real exponents be defined as well. The exponentiation is monotone, which means that for any $x,y\ne\mathbb{0}$, the inequality $x\leq y$ yields the inequality $x^{r}\geq y^{r}$ if $r\leq0$, and $x^{r}\leq y^{r}$ if $r>0$.

{In the algebraic expressions below, the exponents are thought of in the sense of tropical algebra unless otherwise specified. We omit the multiplication symbol $\otimes$ to save writing. }

Examples of the tropical semifields include real semifields that are defined on the sets $\mathbb{R}\cup\{-\infty\}$ and $\mathbb{R}_{+}\cup\{0\}$, where $\mathbb{R}_{+}=\{x\in\mathbb{R}|\ x>0\}$, and given by
\begin{equation*}
\mathbb{R}_{\max,+}
=
(\mathbb{R}\cup\{-\infty\},-\infty,0,\max,+),
\qquad
\mathbb{R}_{\max}
=
(\mathbb{R}_{+}\cup\{0\},0,1,\max,\times).
\end{equation*}

In the semifield $\mathbb{R}_{\max,+}$, which is also known as max-plus algebra, addition $\oplus$ is defined as operation $\max$, and multiplication $\otimes$ as arithmetic addition $+$. The zero $\mathbb{0}$ is set to the number $-\infty$, and the one $\mathbb{1}$ to the arithmetic zero $0$. The inverse $x^{-1}$ is equal to the opposite number $-x$, and the power $x^{y}$ to the product $x\times y$ in standard arithmetic.

{We illustrate the operations in $\mathbb{R}_{\max,+}$ with the next numerical examples (where the symbols $+$ and $\times$ indicate the usual arithmetic addition and multiplication).
\begin{trivlist}
\item
Addition:
\begin{equation*}
\begin{aligned}
&
1\oplus1=1
&&&
(\max(1,1)=1),
\\
&
0\oplus(-3)=0
&&&
(\max(0,-3)=0),
\\
&
2\oplus\mathbb{0}=2
&&&
(\max(2,-\infty)=2).
\end{aligned}
\end{equation*}
\item
Multiplication:
\begin{equation*}
\begin{aligned}
&
1\otimes1=2
&&&
(1+1=2),
\\
&
2\otimes(-3)=-1
&&&
(2+(-3)=-1),
\\
&
1\otimes\mathbb{0}=\mathbb{0}
&&&
(1+(-\infty)=-\infty).
\end{aligned}
\end{equation*}
\item
Exponentiation:
\begin{equation*}
\begin{aligned}
&
1^2=2
&&&
(1\times2=2),
\\
&
(-2)^{1/3}=-2/3
&&&
((-2)\times(1/3)=-2/3),
\\
&
1^{-1}=-1
&&&
(1\times(-1)=-1).
\end{aligned}
\end{equation*}
\end{trivlist}
}

The semifield $\mathbb{R}_{\max}$ (max-algebra) has the operations $\oplus=\max$ and $\otimes=\times$, and the neutral elements $\mathbb{0}=0$ and $\mathbb{1}=1$. The inversion and exponentiation are defined as usual. The order induced by addition in both semifields $\mathbb{R}_{\max,+}$ and $\mathbb{R}_{\max}$ coincides with the natural linear order on $\mathbb{R}$.  Examples are straightforward and hence are omitted.

{Finally, we note that the above semifields are obviously isomorphic to each other by the mappings $\exp:\mathbb{R}_{\max,+}\to\mathbb{R}_{\max}$ and $\log:\mathbb{R}_{\max}\to\mathbb{R}_{\max,+}$.}

\subsection{Distributive Properties}

On the basis of the extremal property of addition, a maximum operation is given by $\max(x,y)=x\oplus y$ for any $x,y\in\mathbb{X}$. A dual minimum operation can then be defined in the following form: $\min(x,y)=(x^{-1}\oplus y^{-1})^{-1}$ if $x,y\ne\mathbb{0}$, and $\min(x,y)=\mathbb{0}$ otherwise.

We now present distributive properties for addition $\oplus$ ($\max$) and $\min$, which are used in the solution of approximation problems below. First, assume that there are scalars $x_{ij}\in\mathbb{X}$ for all $i=1,\ldots,M$ and $j=1,\ldots,N$, where $M$ and $N$ are positive integers. Let $\{I_{1},\ldots,I_{N}\}$ denote a partition that divides the set of naturals $\{1,\ldots,M\}$ into $N$ parts. {The following identity extends the usual distributivity of the operation $\max$ over $\min$ to provide a distributive property of addition $\oplus$ over $\min$ (see, e.g., \cite{Krivulin2024Solution}):}
\begin{equation}
\bigoplus_{i=1}^{M}
\min_{1\leq j\leq N}
x_{ij}
=
\min_{\{I_{1},\ldots,I_{N}\}}
\bigoplus_{j=1}^{N}
\bigoplus_{i\in I_{j}}
x_{ij},
\label{E-plusmin}
\end{equation}
where the minimum on the right-hand side is taken over all partitions $\{I_{1},\ldots,I_{N}\}$.

Suppose that given some functions $f_{j}:\mathbb{X}\rightarrow\mathbb{X}$ for all $j=1,\ldots,N$, we consider a new function defined as the sum $f_{1}(x_{1})\oplus\cdots\oplus f_{N}(x_{N})$. The new function is of the form of max-separable functions, which were investigated in \cite{Zimmermann1984Maxseparable,Zimmermann2003Disjunctive,Tharwat2010One}. Because of the separability, the minimization of this function over all $x_{1},\dots,x_{N}$ reduces to the evaluation of the minimum of each function $f_{j}(x_{j})$, which leads to the following distributive identity: 
\begin{equation}
\min_{x_{1},\ldots,x_{N}}
\bigoplus_{j=1}^{N}
f_{j}(x_{j})
=
\bigoplus_{j=1}^{N}
\min_{x_{j}}
f_{j}(x_{j}).
\label{E-minplus}
\end{equation}

\subsection{Algebra of Matrices and Vectors}

The matrices and vectors with entries in $\mathbb{X}$ are defined in the usual way. The set of matrices of $M$ rows and $N$ columns is denoted by $\mathbb{X}^{M\times N}$. The set of column vectors of $N$ elements is denoted by $\mathbb{X}^{N}$. A matrix (vector) that has all entries equal to $\mathbb{0}$ is the zero matrix (vector) denoted $\bm{0}$. A matrix without zero rows (columns) is called row-regular (column-regular). A matrix is referred to as regular if it is both row- and column-regular.

The operations on matrices (vectors) follow the standard entrywise rules, where the arithmetic addition and multiplication are replaced by the operations $\oplus$ and $\otimes$. In what follows, matrix (vector) multiplication is understood only in the sense of tropical algebra, and the matrix (vector) multiplication sign $\otimes$ is omitted. Specifically, for conforming matrices $\bm{A}=(a_{ij})$, $\bm{B}=(b_{ij})$, $\bm{C}=(c_{ij})$, and scalar $x$, the matrix operations are given by
\begin{equation*}
(\bm{A}\oplus\bm{B})_{ij}
=
a_{ij}\oplus b_{ij},
\qquad
(\bm{A}\bm{C})_{ij}
=
\bigoplus_{k}a_{ik}c_{kj},
\qquad
(x\bm{A})_{ij}
=
xa_{ij}.
\end{equation*}

The monotonicity properties of the scalar operations $\oplus$ and $\otimes$ extend to the matrix (vector) operations where the inequalities are understood entrywise.

Let $\bm{x}=(x_{i})$ be a nonzero column vector. The multiplicative conjugate transpose of $\bm{x}$ is a row vector $\bm{x}^{-}=(x_{i}^{-})$ that has the entries $x_{i}^{-}=x_{i}^{-1}$ if $x_{i}\ne\mathbb{0}$, and $x_{i}^{-}=\mathbb{0}$ otherwise.

The support of a vector $\bm{x}=(x_{i})$ in $\mathbb{X}^{N}$ is defined as $\mathop\mathrm{supp}(\bm{x})=\{i|\ x_{i}\ne\mathbb{0},\ 1\leq i\leq N\}$. For nonzero vectors $\bm{x}=(x_{i})$ and $\bm{y}=(y_{i})$ in $\mathbb{X}^{N}$, such that $\mathop\mathrm{supp}(\bm{x})=\mathop\mathrm{supp}(\bm{y})$, the distance between the vectors is given by the function
\begin{equation}
d(\bm{x},\bm{y})
=
\bigoplus_{i\in\mathop\mathrm{supp}(\bm{x})}\left(y_{i}^{-1}x_{i}\oplus x_{i}^{-1}y_{i}\right)
=
\bm{y}^{-}\bm{x}
\oplus
\bm{x}^{-}\bm{y}.
\label{E-dxy}
\end{equation}

If $\mathop\mathrm{supp}(\bm{x})\ne\mathop\mathrm{supp}(\bm{y})$, then we put $d(\bm{x},\bm{y})=\infty$, where $\infty$ denotes an undefined value greater than any element in $\mathbb{X}$. In the case that $\bm{x}=\bm{y}=\bm{0}$, we assume that $d(\bm{x},\bm{y})=\mathbb{1}$.

In the framework of max-plus algebra $\mathbb{R}_{\max,+}$ where $\mathbb{1}=0$, the function $d$ coincides for all $\bm{x},\bm{y}\in\mathbb{R}^{N}$ with the Chebyshev metric, which is given in standard notation by
\begin{equation*}
d_{\infty}(\bm{x},\bm{y})
=
\max_{1\leq i\leq N}|x_{i}-y_{i}|
=
\max_{1\leq i\leq N}\max(x_{i}-y_{i},y_{i}-x_{i}).
\end{equation*}

In the case of max-algebra $\mathbb{R}_{\max}$, the function $d$ can be considered as a generalized metric with values in the interval $[1,\infty)\subset\mathbb{R}_{+}$. Due to the isomorphism between $\mathbb{R}_{\max}$ and $\mathbb{R}_{\max,+}$, this function can be transformed into a metric $d^{\prime}(\bm{x},\bm{y})=\log d(\bm{x},\bm{y})$.

\subsection{Best Approximate Solution to Vector Equations}

Given a matrix $\bm{A}\in\mathbb{X}^{M\times N}$ and vector $\bm{b}\in\mathbb{X}^{M}$, consider the problem of finding vectors $\bm{x}\in\mathbb{X}^{N}$ that satisfy the equation with the unknown on one side in the form
\begin{equation}
\bm{A}\bm{x}
=
\bm{b},
\label{E-Axeqb}
\end{equation}
where the left-hand matrix-vector multiplication is performed in the tropical algebra sense.

Solutions to this one-sided equation are known in various forms depending on the assumptions and technique used. Since Equation \eqref{E-Axeqb} may have no solution, one can consider finding a best approximate solution in the sense of the metric $d$. A vector $\bm{x}_{\ast}$ is a best approximate solution of the equation if for all vectors $\bm{x}$ the following inequality holds:
\begin{equation*}
d(\bm{A}\bm{x}_{\ast},\bm{b})
\leq
d(\bm{A}\bm{x},\bm{b}).
\end{equation*}

{The next approximate solution is given in \cite{Krivulin2009Onsolution,Krivulin2012Solution}.}

\begin{theorem}
\label{T-Axeqb}
Let $\bm{A}$ be a regular matrix and $\bm{b}$ a regular vector. Define the scalar $\Delta=(\bm{A}(\bm{b}^{-}\bm{A})^{-})^{-}\bm{b}$. Then, the following statements hold:
\begin{enumerate}
\item
The best approximation error for equation \eqref{E-Axeqb} is equal to
\begin{equation*}
d(\bm{A}\bm{x}_{\ast},\bm{b})
=
\sqrt{\Delta};
\end{equation*}
\item
The best approximate solution of \eqref{E-Axeqb} is given by
\begin{equation*}
\bm{x}_{\ast}
=
\sqrt{\Delta}(\bm{b}^{-}\bm{A})^{-}.
\end{equation*}
\item
If $\Delta=\mathbb{1}$, then the equation has exact solutions, including the maximal solution $\bm{x}_{\ast}=(\bm{b}^{-}\bm{A})^{-}$.
\end{enumerate}
\end{theorem}

The computational complexity of the solution is at most $O(MN)$. 

Suppose given matrices $\bm{A}\in\mathbb{X}^{M\times N}$ and $\bm{B}\in\mathbb{X}^{M\times L}$, we need to find regular vectors $\bm{x}\in\mathbb{X}^{N}$ and $\bm{y}\in\mathbb{X}^{L}$ that solve the two-sided equation
\begin{equation}
\bm{A}\bm{x}
=
\bm{B}\bm{y}.
\label{E-AxeqBy}
\end{equation}

An alternating algorithm to obtain a best approximate solution to equation \eqref{E-AxeqBy} in the case of regular matrices is proposed and its convergence is examined in \cite{Krivulin2023Solution}, where references to other existing solutions of the problem are also given. The algorithm is described as an iterative computational procedure in the following form (Algorithm~\ref{A-AxeqBy}).

\begin{algorithm}[H]
\caption{{Approximate}
 solution of two-sided equation $\bm{A}\bm{x}=\bm{B}\bm{y}$.}
\label{A-AxeqBy}
\begin{enumerate}
\item
Input regular matrices $\bm{A},\bm{B}$ and regular vector $\bm{x}_{0}$; put $k=0$.
\item\label{Loop-AxeqBy}
Calculate the squared approximation error and vector
\begin{equation*}
\Delta_{k}
=
(\bm{B}((\bm{A}\bm{x}_{k})^{-}\bm{B})^{-})^{-}\bm{A}\bm{x}_{k},
\qquad
\bm{y}_{k+1}
=
\sqrt{\Delta_{k}}((\bm{A}\bm{x}_{k})^{-}\bm{B})^{-}.
\end{equation*}
\item
If $\Delta_{k}=\mathbb{1}$ or $\bm{y}_{k+1}=\bm{y}_{j}$ for some $j<k$, then set
\begin{equation*}
\Delta_{\ast}
=
\Delta_{k},
\qquad
\bm{x}_{\ast}
=
\bm{x}_{k},
\qquad
\bm{y}_{\ast}
=
\bm{y}_{k+1},
\end{equation*}
and stop; otherwise set $k$ to $k+1$.
\item
Calculate the squared approximation error and vector
\begin{equation*}
\Delta_{k}
=
(\bm{A}((\bm{B}\bm{y}_{k})^{-}\bm{A})^{-})^{-}\bm{B}\bm{y}_{k},
\qquad
\bm{x}_{k+1}
=
\sqrt{\Delta_{k}}((\bm{B}\bm{y}_{k})^{-}\bm{A})^{-}.
\end{equation*}
\item
If $\Delta_{k}=\mathbb{1}$ or $\bm{x}_{k+1}=\bm{x}_{j}$ for some $j<k$, then set
\begin{equation*}
\Delta_{\ast}
=
\Delta_{k},
\qquad
\bm{x}_{\ast}
=
\bm{x}_{k+1},
\qquad
\bm{y}_{\ast}
=
\bm{y}_{k},
\end{equation*}
and stop; otherwise set $k$ to $k+1$.
\item
Go to step \ref{Loop-AxeqBy}.
\end{enumerate}
\end{algorithm}
\vspace{+9pt}

As the algorithm output, one obtains the approximation error $\Delta_{\ast}$ and approximate solutions $\bm{x}_{\ast}$ and $\bm{y}_{\ast}$. If $\Delta_{\ast}=\mathbb{1}$, then the obtained vectors $\bm{x}_{\ast}$ and $\bm{y}_{\ast}$ are an exact solution.

We observe that calculation of the squared error and solution vector at each iteration requires no more than $O(M\max(N,L))$ operations.

\section{Tropical Polynomials and Rational Functions}
\label{S-TPRF}

Tropical polynomials are defined as analogues of the polynomials in the usual setting, where the ordinary addition, multiplication, and exponentiation are replaced by their tropical counterparts. In the approximation problems below, we use tropical analogues of generalized Puiseux polynomials in which the exponents can take any real values. 

\subsection{Puiseux Polynomial Functions}

We consider a (generalized) Puiseux polynomial over $\mathbb{X}$ with $N$ monomials in one variable $x$, which is given by
\begin{equation}
P(x)
=
\bigoplus_{j=1}^{N}\theta_{j}x^{p_{j}}
=
\theta_{1}x^{p_{1}}\oplus\cdots\oplus\theta_{N}x^{p_{N}},
\qquad
x\ne\mathbb{0},
\label{E-Px}
\end{equation}
where $p_{1},\ldots,p_{N}\in\mathbb{R}$ are exponents and $\theta_{1},\ldots,\theta_{N}\in\mathbb{X}$ are nonzero coefficients.

If interpreted in terms of max-plus algebra $\mathbb{R}_{\max,+}$, the polynomial can be expressed using standard arithmetic operations in the form
\begin{equation*}
P(x)
=
\max_{1\leq j\leq N}(p_{j}\times x+\theta_{j}),
\qquad
x\in\mathbb{R},
\end{equation*}
with $p_{1},\ldots,p_{N}\in\mathbb{R}$ and $\theta_{1},\ldots,\theta_{N}\in\mathbb{R}$, and defines a piecewise linear convex function. 

In the context of max-algebra $\mathbb{R}_{\max}$, the polynomial is represented using ordinary operations as a spline function
\begin{equation*}
P(x)
=
\max_{1\leq j\leq N}(\theta_{j}\times x^{p_{j}}),
\qquad
x\in\mathbb{R}_{+},
\end{equation*}
where $p_{1},\ldots,p_{N}\in\mathbb{R}_{+}$ and $\theta_{1},\ldots,\theta_{N}\in\mathbb{R}_{+}$, which may not be convex.

We observe that the isomorphism between the semifields $\mathbb{R}_{\max,+}$ and $\mathbb{R}_{\max}$ provides a one-to-one correspondence between polynomials given in terms of these semifields. Specifically, taking a logarithm to a base greater than one, which is monotone increasing, transforms a max-algebra polynomial $P(x)$ into a max-plus algebra polynomial as
\begin{equation*}
P^{\prime}(x)
=
\log P(x)
=
\max_{1\leq j\leq N}(\log\theta_{j}+p_{j}\times\log x)
=
\max_{1\leq j\leq N}(p_{j}\times x^{\prime}+\theta_{j}^{\prime})
\end{equation*}
with new indeterminate $x^{\prime}=\log x$ and coefficients $\theta_{j}^{\prime}=\log\theta_{j}$ for all $j=1,\ldots,N$.

As a result, the calculation of polynomials in the context of max-algebra can be directly reduced to that in terms of max-plus algebra and vice versa. 

We conclude this subsection with a polynomial optimization problem that is involved in the approximation procedure presented below. Consider a minimization problem of a polynomial function \eqref{E-Px} formulated in the form
\begin{equation}
\min_{x>\mathbb{0}}\ 
\bigoplus_{j=1}^{N}
\theta_{j}x^{p_{j}},
\label{P-minx_sumjthetajxpj}
\end{equation}
where the polynomial is assumed to have exponents $p_{j}$ of different signs.

{To solve the problem, one can apply the next result obtained in \cite{Krivulin2021Algebraic} (see also \cite{Krivulin2024Solution}).}
\begin{lemma}
\label{L-minx_sumjthetajxpj}
The minimum in problem \eqref{P-minx_sumjthetajxpj} is given by
\begin{equation*}
\mu
=
\bigoplus_{\substack{1\leq j,k\leq N\\p_{j}<0,\ p_{k}>0}}
\theta_{j}^{-\frac{p_{k}}{p_{j}-p_{k}}}\theta_{k}^{\frac{p_{j}}{p_{j}-p_{k}}}
\oplus
\bigoplus_{\substack{1\leq j\leq N\\p_{j}=0}}
\theta_{j},
\end{equation*}
and all solutions of \eqref{P-minx_sumjthetajxpj} satisfy the condition
\begin{equation*}
\bigoplus_{\substack{1\leq j\leq N\\p_{j}<0}}
\mu^{1/p_{j}}
\theta_{j}^{-1/p_{j}}
\leq
x
\leq
\min_{\substack{1\leq j\leq N\\p_{j}>0}}
\mu^{1/p_{j}}
\theta_{j}^{-1/p_{j}}.
\end{equation*}
\end{lemma}

This solution has computational complexity that grows not faster than $O(N^{2})$.

\subsection{Puiseux Rational Functions}

Suppose two polynomials of $N$ and $L$ monomials are given in the form
\begin{equation*}
P(x)
=
\bigoplus_{j=1}^{N}\theta_{j}x^{p_{j}},
\qquad
Q(x)
=
\bigoplus_{j=1}^{L}\sigma_{j}x^{q_{j}}.
\end{equation*}

Taking these polynomials as a numerator and denominator, we define a (generalized) Puiseux rational function of one variable $x$ to be 
\begin{equation}
R(x)
=
\frac{P(x)}{Q(x)}
=
\frac{\theta_{1}x^{p_{1}}\oplus\cdots\oplus\theta_{N}x^{p_{N}}}{\sigma_{1}x^{q_{1}}\oplus\cdots\oplus\sigma_{L}x^{q_{L}}},
\qquad
x\ne\mathbb{0}.
\label{E-Rx}
\end{equation}

In the setting of max-plus algebra $\mathbb{R}_{\max,+}$, the rational function appears to be a difference of two piecewise linear convex functions, whereas in max-algebra $\mathbb{R}_{\max}$, this function becomes a ratio between two spline functions.  

The class of real difference-of-convex (DC) functions is known to be very rich \cite{Hartman1959Onfunctions,Tuy2016Convex}. Specifically, any continuous function can be approximated by a DC function, which makes the approximation by Puiseux rational functions theoretically and practically sound.

We note that, as in the case of polynomials, any tropical rational function given in terms of max-algebra can be directly transformed into one in the max-plus algebra setting.

\section{Tropical Best Discrete Approximation Problems}
\label{TBDAP}

We consider  the best discrete approximation problems that are formulated and solved in the framework of tropical algebra. Suppose that the sample data of $M$ values $y_{1},\ldots,y_{M}$ at points $x_{1},\ldots,x_{M}$ are given for an unknown function $f:\mathbb{X}\rightarrow\mathbb{X}$. The discrete approximation of the function $f(x)$ is to fit the data to a parametric function $F_{\bm{\theta}}:\mathbb{X}\rightarrow\mathbb{X}$ by finding a vector of parameters $\bm{\theta}$ that achieves minimal difference between both sides of the equations 
\begin{equation*}
F_{\bm{\theta}}(x_{i})
=
y_{i},
\qquad
i=1,\ldots,M.
\end{equation*}

We use the vector notation $\bm{x}=(x_{i})$, $\bm{y}=(y_{i})$ and $\bm{F}_{\bm{\theta}}(\bm{x})=(F_{\bm{\theta}}(x_{i}))$, and then arrive at a problem of best approximation in the sense of the distance function $d$ in the form
\begin{equation*}
\min_{\bm{\theta}}\ 
d(\bm{F}_{\bm{\theta}}(\bm{x}),\bm{y}).
\end{equation*}

An approximate solution is found as a minimizer
\begin{equation*}
\bm{\theta}_{\ast}
=
\arg\min_{\bm{\theta}}d(\bm{F}_{\bm{\theta}}(\bm{x}),\bm{y}).
\end{equation*}

Below we outline problems of best discrete approximation by Puiseux polynomials and rational functions and discuss solution approaches under various assumptions.

\subsection{Approximation by Polynomial Functions}
Consider the approximation of an unknown function $f(x)$ by a polynomial function $P(x)$ given by \eqref{E-Px}. We formulate the problem to find both exponents $p_{1},\ldots,p_{N}$ and coefficients $\theta_{1},\ldots,\theta_{N}$ that provide the best agreement between two sides of the equations
\begin{equation}
\theta_{1}x_{i}^{p_{1}}\oplus\cdots\oplus\theta_{N}x_{i}^{p_{N}}
=
y_{i},
\qquad
i=1,\ldots,M.
\label{E-theta1xip1-yi}
\end{equation}

To solve the problem, we first introduce the vectors and matrix
\begin{equation*}
\bm{y}
=
\left(
\begin{array}{c}
y_{1}
\\
\vdots
\\
y_{M}
\end{array}
\right),
\qquad
\bm{p}
=
\left(
\begin{array}{c}
p_{1}
\\
\vdots
\\
p_{N}
\end{array}
\right),
\qquad
\bm{\theta}
=
\left(
\begin{array}{c}
\theta_{1}
\\
\vdots
\\
\theta_{N}
\end{array}
\right),
\qquad
\bm{X}(\bm{p})
=
\left(
\begin{array}{ccc}
x_{1}^{p_{1}} & \ldots & x_{1}^{p_{N}}
\\
\vdots & & \vdots
\\
x_{M}^{p_{1}} & \ldots & x_{M}^{p_{N}}
\end{array}
\right).
\end{equation*}

With this notation, equations \eqref{E-theta1xip1-yi} can be rewritten in vector form as
\begin{equation}
\bm{X}(\bm{p})\bm{\theta}
=
\bm{y}.
\label{E-Xpthetaeqy}
\end{equation}

We formulate the problem as the best approximation of the vector equation at \eqref{E-Xpthetaeqy} in the sense of the distance function $d$ to find
\begin{equation*}
(\bm{p}_{\ast},\bm{\theta}_{\ast})
=
\arg\min_{\bm{p},\bm{\theta}}d(\bm{X}(\bm{p})\bm{\theta},\bm{y}).
\end{equation*}

If the vector of exponents $\bm{p}$ is fixed in advance, then the entries in the matrix $\bm{X}(\bm{p})$ are completely defined. As a result, a solution of the approximation problem can be obtained by a direct application of Theorem~\ref{T-Axeqb}. The solution takes no more than $O(MN)$ operations to calculate the squared approximation error and the vector of coefficients given by
\begin{equation}
\Delta_{\ast}
=
(\bm{X}(\bm{p})(\bm{y}^{-}\bm{X}(\bm{p}))^{-})^{-}\bm{y},
\qquad
\bm{\theta}_{\ast}
=
\sqrt{\Delta_{\ast}}(\bm{y}^{-}\bm{X}(\bm{p}))^{-}.
\label{E-Deltaast-thetaast}
\end{equation}

In the general case, when both vectors of exponents $\bm{p}$ and coefficients $\bm{\theta}$ are unknown, one can apply the procedure proposed in \cite{Krivulin2023Algebraic}, which combines random search over a set of vectors of exponents with application of Theorem~\ref{T-Axeqb} to obtain a vector of coefficients. The procedure performs a number of iterations, each consisting of sampling a new random vector of exponents along with evaluating the error and the vector of coefficients according to \eqref{E-Deltaast-thetaast}. However, as numerical experience shows, the number of iterations required to achieve a fairly accurate solution is quite large even for problems of moderate dimensions  and becomes unacceptable as the number of monomials $N$ increases. 

To overcome this difficulty, a more efficient solution approach is developed in \cite{Krivulin2024Solution}. The approach addresses approximation problems set in terms of max-plus algebra $\mathbb{R}_{\max,+}$ but can easily be extended to other tropical semifields isomorphic to $\mathbb{R}_{\max,+}$, including max-algebra $\mathbb{R}_{\max}$. The solution starts with an application of Theorem~\ref{T-Axeqb} to vector Equation \eqref{E-Xpthetaeqy} where $\bm{X}(\bm{p})$ is considered a matrix parameterized by the vector $\bm{p}$. As a result, we obtain an approximate solution represented in parametric form as $\bm{\theta}(\bm{p})=\sqrt{\delta(\bm{p})}(\bm{y}^{-}\bm{X}(\bm{p}))^{-}$ where $\delta(\bm{p})=(\bm{X}(\bm{p})(\bm{y}^{-}\bm{X}(\bm{p}))^{-})^{-}\bm{y}$ is the squared approximation error.

To complete the solution, we need to minimize the error $\delta(\bm{p})$ with respect to $\bm{p}$, which yields the optimal vector of exponents
\begin{equation}
\bm{p}_{\ast}
=
\arg\min_{\bm{p}}(\bm{X}(\bm{p})(\bm{y}^{-}\bm{X}(\bm{p}))^{-})^{-}\bm{y}.
\label{E-past}
\end{equation}

After calculation of $\bm{p}_{\ast}$, one can evaluate the minimal squared approximation error $\Delta_{\ast}=\delta(\bm{p}_{\ast})$ and its related vector of coefficients $\bm{\theta}_{\ast}=\sqrt{\Delta_{\ast}}(\bm{y}^{-}\bm{X}(\bm{p}_{\ast}))^{-}$.

We find the minimum at \eqref{E-past} by transformation into a combinatorial optimization problem and application of a solution procedure based on the agglomerative clustering technique. This procedure requires computational complexity of order no more than $O(\max(M^{4},(M-N)^{3}M^{3}))$. We describe the solution in more detail in the next section.

\subsection{Approximation by Rational Functions}

Consider an approximation problem with the same sample data as above, where the unknown function $f(x)$ needs to be approximated by a rational function $R(x)$ given by \eqref{E-Rx}. We now seek two sets of exponents $p_{1},\ldots,p_{N}$ and $q_{1},\ldots,q_{L}$ and two sets of parameters $\theta_{1},\ldots,\theta_{N}$ and $\sigma_{1},\ldots,\sigma_{L}$ to obtain the best approximation of the equations 
\begin{equation}
\frac{\theta_{1}x_{i}^{p_{1}}\oplus\cdots\oplus\theta_{N}x_{i}^{p_{N}}}{\sigma_{1}x_{i}^{q_{1}}\oplus\cdots\oplus\sigma_{L}x_{i}^{q_{L}}}
=
y_{i},
\qquad
i=1,\ldots,M.
\label{E-theta1xip1_sigma1xiq1_yi}
\end{equation}

To represent Equation \eqref{E-theta1xip1_sigma1xiq1_yi} in vector form, we first rewrite these equations as
\begin{equation}
\theta_{1}x_{i}^{p_{1}}\oplus\cdots\oplus\theta_{N}x_{i}^{p_{N}}
=
y_{i}(\sigma_{1}x_{i}^{q_{1}}\oplus\cdots\oplus\sigma_{L}x_{i}^{q_{L}}),
\qquad
i=1,\ldots,M.
\label{E-theta1xip1-yisigma1xiq1}
\end{equation}

In addition to the notations $\bm{p}$, $\bm{\theta}$ and $\bm{X}(\bm{p})$ introduced above, we define the following vectors and matrices:
\begin{equation*}
\bm{q}
=
\left(
\begin{array}{c}
q_{1}
\\
\vdots
\\
q_{L}
\end{array}
\right),
\quad
\bm{\sigma}
=
\left(
\begin{array}{c}
\sigma_{1}
\\
\vdots
\\
\sigma_{L}
\end{array}
\right),
\quad
\bm{Y}
=
\left(
\begin{array}{ccc}
y_{1} &  & \mathbb{0}
\\
& \ddots &
\\
\mathbb{0} & & y_{M}
\end{array}
\right),
\quad
\bm{Z}(\bm{q})
=
\left(
\begin{array}{ccc}
x_{1}^{q_{1}} & \ldots & x_{1}^{q_{L}}
\\
\vdots & & \vdots
\\
x_{M}^{q_{1}} & \ldots & x_{M}^{q_{L}}
\end{array}
\right).
\end{equation*}

In vector notation, Equation \eqref{E-theta1xip1-yisigma1xiq1} takes the form
\begin{equation}
\bm{X}(\bm{p})\bm{\theta}
=
\bm{Y}\bm{Z}(\bm{q})\bm{\sigma}.
\label{E-Xptheta_YZqsigma}
\end{equation}

One can verify that approximation error for Equation \eqref{E-Xptheta_YZqsigma} written in terms of the distance function $d$ coincides with the error of the system at \eqref{E-theta1xip1_sigma1xiq1_yi} (see \cite{Krivulin2023Algebraic}).

Finally, we formulate a best approximation problem  to find
\begin{equation}
(\bm{p}_{\ast},\bm{\theta}_{\ast},\bm{q}_{\ast},\bm{\sigma}_{\ast})
=
\arg\min_{\bm{p},\bm{\theta},\bm{q},\bm{\sigma}}d(\bm{X}(\bm{p})\bm{\theta},\bm{Y}\bm{Z}(\bm{q})\bm{\sigma}).
\end{equation}

To discuss solutions to problem \eqref{E-Xptheta_YZqsigma}, first assume that both vectors of exponents $\bm{p}$ and $\bm{q}$ are set to some predetermined values, which makes both parameterized matrices $\bm{X}(\bm{p})$ and $\bm{Z}(\bm{q})$ become fixed. In this case, the problem is solved by applying Algorithm~\ref{A-AxeqBy} , which produces the squared approximation error $\Delta_{\ast}$ and two vectors of coefficients $\bm{\theta}_{\ast}$ and $\bm{\sigma}_{\ast}$.
 
Suppose now that the vectors of exponents $\bm{p}$ and $\bm{q}$ are unknown and must be evaluated along with the unknown vectors of coefficients $\bm{\theta}$ and $\bm{\sigma}$. One of the approaches to handle this approximation problem is to apply a random search procedure described in~\cite{Krivulin2023Algebraic}. The procedure combines random sampling over all possible vectors of exponents with alternate iterations according to Algorithm~\ref{A-AxeqBy} to obtain corresponding vectors of coefficients and then calculate the squared approximation error. Those vectors of exponents and coefficients that  lead to the minimal error $\Delta_{\ast}$ over all samples are taken as solution vectors $\bm{p}_{\ast}$, $\bm{q}_{\ast}$, $\bm{\theta}_{\ast}$ and $\bm{\sigma}_{\ast}$ that determine an approximate rational function.

We observe that this random search procedure is based on repeated executions of Algorithm~\ref{A-AxeqBy} and thus can be much more expensive than the random search procedure for the polynomial approximation described above.

As an alternative to the computationally inefficient random search technique, we propose a procedure that combines alternate calculations according to Algorithm~\ref{A-AxeqBy} with simultaneous evaluation of exponents and coefficients as described in \cite{Krivulin2024Solution}. The procedure is developed in the context of max-plus algebra  but can be easily extended to solve approximation problems in terms of other isomorphic semifields. Further details on the proposed procedure are given in the subsequent sections.

\section{Approximation by Polynomial Functions}
\label{S-APF}

We start with the solution proposed in \cite{Krivulin2024Solution} in the framework of the semifield $\mathbb{R}_{\max,+}$ (max-plus algebra) for polynomial approximation problems, where both exponents and coefficients of monomials in the approximating polynomial are unknown and thus need to be evaluated. This solution offers an alternating computational scheme of solving a sequence of polynomial approximations derived from the initial rational approximation problem and plays a pivotal role in the new technique of rational approximation.

The solution of the polynomial approximation problem reduces to solving Equation \eqref{E-Xpthetaeqy} for both unknown vectors $\bm{p}$ and $\bm{\theta}$. We implement the procedure developed in \cite{Krivulin2024Solution}, which includes two main stages. The procedure starts with an application of Theorem~\ref{T-Axeqb} as the first stage to produce a direct representation of the squared approximation error and the vector of coefficients, both parameterized by the vector of exponents in the form
\begin{equation}
\delta(\bm{p})
=
(\bm{X}(\bm{p})(\bm{y}^{-}\bm{X}(\bm{p}))^{-})^{-}\bm{y},
\qquad
\bm{\theta}(\bm{p})
=
\sqrt{\delta(\bm{p})}(\bm{y}^{-}\bm{X}(\bm{p}))^{-}.
\label{E-deltap-thetap}
\end{equation}

At the second stage, the procedure minimizes the error to find the optimal vector of exponents $\bm{p}_{\ast}$ by solving the minimization problem
\begin{equation}
\min_{\bm{p}}\ 
\delta(\bm{p}).
\label{P-minpdeltap}
\end{equation}

Substitution of the obtained vector $\bm{p}_{\ast}$ into \eqref{E-deltap-thetap} yields $\Delta_{\ast}=\delta(\bm{p}_{\ast})$ and $\bm{\theta}_{\ast}=\bm{\theta}(\bm{p}_{\ast})$ and thus completes the solution. In what follows, we describe the solution in more detail. Specifically, we outline the transformation technique of the error function and include the algorithm of minimizing the error from \cite{Krivulin2024Solution} for the sake of completeness.

\subsection{Minimization of Approximation Error}

To solve the minimization problem at \eqref{P-minpdeltap}, we exploit the representation of $\min$ through $\oplus$ and rewrite the objective function as follows:
\begin{equation*}
\delta(\bm{p})
=
(\bm{X}(\bm{p})(\bm{y}^{-}\bm{X}(\bm{p}))^{-})^{-}\bm{y}
=
\bigoplus_{i=1}^{M}
\left(
\bigoplus_{j=1}^{N}
(\varphi_{i}(p_{j}))^{-1}
\right)^{-1}
=
\bigoplus_{i=1}^{M}
\min_{1\leq j\leq N}
\varphi_{i}(p_{j}),
\end{equation*}
where the symbols $\varphi_{i}$ stand for the functions
\begin{equation*}
\varphi_{i}(p)
=
y_{i}x_{i}^{-p}(y_{1}^{-1}x_{1}^{p}\oplus\cdots\oplus y_{M}^{-1}x_{M}^{p}),
\qquad
i=1,\ldots,M.
\end{equation*}

Applying identity \eqref{E-plusmin} puts the objective function into the form
\begin{equation*}
\delta(\bm{p})
=
\min_{\{I_{1},\ldots,I_{N}\}} 
\bigoplus_{j=1}^{N}
\bigoplus_{i\in I_{j}}
\varphi_{i}(p_{j}),
\end{equation*}
where the minimum is over all partitions $\{I_{1},\ldots,I_{N}\}$ of the set $\{1,\ldots,M\}$ into $N$ parts.

After the substitution of the obtained representation for the function $\delta(\bm{p})$ and the change of the order of minimization, the problem at \eqref{P-minpdeltap} becomes
\begin{equation*}
\min_{\{I_{1},\ldots,I_{N}\}}\ 
\min_{p_{1},\ldots,p_{N}}
\bigoplus_{j=1}^{N}
\bigoplus_{i\in I_{j}}
\varphi_{i}(p_{j}).
\end{equation*}

The objective function in the last problem is the minimum of a max-separable function, which can be rearranged according to identity \eqref{E-minplus} to rewrite the problem as 
\begin{equation}
\min_{\{I_{1},\ldots,I_{N}\}}\ 
\bigoplus_{j=1}^{N}
\min_{p_{j}}
\bigoplus_{i\in I_{j}}
\varphi_{i}(p_{j}).
\label{P-minI1INminpjvarphiipj}
\end{equation}

The evaluation of the objective function in problem \eqref{P-minI1INminpjvarphiipj} for each partition $\{I_{1},\ldots,I_{N}\}$ involves solving $N$ inner minimization problems to find $p_{1}^{\ast},\ldots,p_{N}^{\ast}$ as their minimizers. The maximum (tropical sum) of the minimum values in the inner problems defines the corresponding value of the objective function in \eqref{P-minI1INminpjvarphiipj}. The set of minimizers that produce the minimum of the objective function $\delta(\bm{p})$ over all partitions presents a vector of optimal exponents $\bm{p}_{\ast}$, whereas the value $\Delta_{\ast}=\delta(\bm{p}_{\ast})$ gives the best squared approximation error.

The solution of each inner minimization problem is obtained in terms of the semifield $\mathbb{R}_{\max,+}$ as follows. In this semifield, the identity $x^{y}=y^{x}$ holds for all $x,y\in\mathbb{R}$ since both sides of the identity correspond to the arithmetic product $x\times y$. Based on this identity, the functions $\varphi_{i}(p)$ transform into polynomials, where $p$ becomes the indeterminate, to write
\begin{equation*}
\varphi_{i}(p)
=
y_{i}p^{-x_{i}}(y_{1}^{-1}p^{x_{1}}\oplus\cdots\oplus y_{M}^{-1}p^{x_{M}}),
\qquad
i=1,\ldots,M.
\end{equation*}

As a result, the inner minimization problems in \eqref{P-minI1INminpjvarphiipj} turn to polynomial optimization problems with respect to $p$ that are solved by application of Lemma~\ref{L-minx_sumjthetajxpj}.

To solve the outer optimization problem at \eqref{P-minI1INminpjvarphiipj}, we apply a computational scheme based on an agglomerative clustering technique, which finds an optimal (near-optimal) partition. The solution starts with the partition of the set $\{1,\ldots,M\}$ into $M$ single-element subsets and then iteratively reduces the number of partitions by merging subsets in pairs. 

{Each subset in a partition defines a polynomial whose minimum value becomes a characteristic of the subset. The maximum of the minimums over all subsets produces the value of the objective function corresponding to the partition. At each step of the scheme, we merge those subsets in the partition for which the merged subset yields a polynomial with the least minimum value over all pairs of subsets. The merging process continues until the number of subsets in a new partition becomes equal to $N$.}

\subsection{Error Minimization Procedure}

We give a formal description of the overall procedure, which uses sample data to evaluate the optimal approximation error and find corresponding exponents in the framework of the semifield $\mathbb{R}_{\max,+}$ in the form of Algorithm~\ref{A-MFdeltap}. 

\vspace{+9pt}

\begin{algorithm}[H]
 
\caption{{Minimization} of error function $\delta(\bm{p})$ in max-plus algebra.}
\label{A-MFdeltap}
\begin{enumerate}
\item
Input $M$ samples $(x_{i},y_{i})$ for $i=1,\ldots,M$, fix positive integer $N$ and put $k=1$.

Define the polynomials \vspace{-4pt}
\begin{equation*}
\varphi_{i}(p)
=
\bigoplus_{j=1}^{M}
y_{j}^{-1}y_{i}
p^{x_{j}-x_{i}},
\qquad
i=1,\ldots,M.
\end{equation*} 

Construct the subsets and form the partition 
\begin{equation*}
I_{1}^{(0)}
=
\{1\},
\ldots,
I_{M}^{(0)}
=
\{M\};
\qquad
\mathcal{P}^{(0)}
=
\{I_{1}^{(0)},\ldots,I_{M}^{(0)}\}.
\end{equation*}
\item\label{L-MFdeltap}
Apply Lemma~\ref{L-minx_sumjthetajxpj} to select a pair of subsets to satisfy the condition \vspace{-3pt}
\begin{equation*}
(U^{(k)},V^{(k)})
=
\arg\min_{\substack{U,V\in\mathcal{P}^{(k-1)}\\ U\ne V}}
\min_{p}
\left(
\bigoplus_{i\in U}
\varphi_{i}(p)
\oplus
\bigoplus_{i\in V}
\varphi_{i}(p)
\right).
\end{equation*}
\item
Merge the subsets $U^{(k)}$ and $V^{(k)}$ to form the partition 
\begin{equation*}
\mathcal{P}^{(k)}
=
(\mathcal{P}^{(k-1)}\setminus\{U^{(k)},V^{(k)}\})\cup\{U^{(k)}\cup V^{(k)}\}.
\end{equation*}
\item
If $M-k=N$, then define the partition
\begin{equation*}
\mathcal{P}_{\ast}
=
\{I_{1}^{\ast},\ldots,I_{N}^{\ast}\}
=
\mathcal{P}^{(M-k)}
=
\{I_{1}^{(M-k)},\ldots,I_{N}^{(M-k)}\};
\end{equation*}
otherwise set $k$ to $k+1$ and go to Step~\ref{L-MFdeltap}.
\item
Apply Lemma~\ref{L-minx_sumjthetajxpj} to find the minimums and corresponding solutions
\begin{gather*}
\delta_{1}^{\ast}
=
\min_{p}
\bigoplus_{i\in I_{1}^{\ast}}
\varphi_{i}(p),
\quad
\ldots,
\quad
\delta_{N}^{\ast}
=
\min_{p}
\bigoplus_{i\in I_{N}^{\ast}}
\varphi_{i}(p);
\\
p_{1}^{\ast}
=
\arg\min_{p}
\bigoplus_{i\in I_{1}^{\ast}}
\varphi_{i}(p),
\quad
\ldots,
\quad
p_{N}^{\ast}
=
\arg\min_{p}
\bigoplus_{i\in I_{N}^{\ast}}
\varphi_{i}(p).
\end{gather*}
\item
Evaluate the minimum of the function $\delta(\bm{p})$ to be\vspace{-4pt}
\begin{equation*}
\Delta_{\ast}
=
\bigoplus_{j=1}^{N}
\delta_{j}^{\ast}
=
\max_{1\leq j\leq N}
\delta_{j}^{\ast}. \vspace{-9pt}
\end{equation*}
\end{enumerate}
\end{algorithm}
\subsection{Polynomial Approximation Procedure}

We apply Algorithm~\ref{A-MFdeltap} as a key ingredient in the solution of the polynomial approximation problems in terms of the semifield $\mathbb{R}_{\max,+}$. The solution procedure is as follows (Algorithm~\ref{A-AP}).

\vspace{+12pt}

\begin{algorithm}[H]
\caption{{Approximation}
 by polynomials in max-plus algebra.}
\label{A-AP}
\begin{enumerate}
\item
Input $M$ samples $(x_{i},y_{i})$ for $i=1,\ldots,M$, and fix positive integer $N$.

Define the vectors and parameterized matrix
\begin{equation*}
\bm{y}
=
\left(
\begin{array}{c}
y_{1}
\\
\vdots
\\
y_{M}
\end{array}
\right),
\qquad
\bm{p}
=
\left(
\begin{array}{c}
p_{1}
\\
\vdots
\\
p_{N}
\end{array}
\right),
\qquad
\bm{X}(\bm{p})
=
\left(
\begin{array}{ccc}
x_{1}^{p_{1}} & \ldots & x_{1}^{p_{N}}
\\
\vdots & & \vdots
\\
x_{M}^{p_{1}} & \ldots & x_{M}^{p_{N}}
\end{array}
\right).
\end{equation*}

Construct the function
\begin{equation*}
\delta(\bm{p})
=
(\bm{X}(\bm{p})(\bm{y}^{-}\bm{X}(\bm{p}))^{-})^{-}\bm{y}.
\end{equation*}
\item
Apply Algorithm~\ref{A-MFdeltap} to find the vector of exponents and squared error
\begin{equation*}
\bm{p}_{\ast}
=
\arg\min_{\bm{p}}
\delta(\bm{p}),
\qquad
\Delta_{\ast}
=
\delta(\bm{p}_{\ast}).
\end{equation*}
\item
Calculate the vector of coefficients
\begin{equation*}
\bm{\theta}_{\ast}
=
\sqrt{\Delta_{\ast}}
(\bm{y}^{-}\bm{X}(\bm{p}_{\ast}))^{-}.\vspace{-6pt}
\end{equation*}
\end{enumerate}
\end{algorithm}

\vspace{+12pt}
 
To conclude this section, we note that the polynomial approximation procedure, including the error minimization algorithm, assumes the max-plus algebra setting of the polynomials involved. However, due to the isomorphism between the semifields $\mathbb{R}_{\max,+}$ and $\mathbb{R}_{\max}$, the approximation problems given in the framework of max-algebra can be readily transformed into their counterparts in terms of max-plus algebra. This makes the procedure readily applicable to solve approximation problems by max-algebra polynomials.

\section{Approximation by Rational Functions}
\label{S-ARF}

In this section, we present a new computational approach to solve rational approximation problems given in the framework of the semifield $\mathbb{R}_{\max,+}$ under the assumption that both exponents and coefficients in the numerator and denominator polynomials are unknown. The approach incorporates solving polynomial approximation problems according to Algorithm~\ref{A-AP} into the alternating procedure based on Algorithm~\ref{A-AxeqBy}.

We consider the approximation problem represented in the form of \eqref{E-Xptheta_YZqsigma}, where the vectors of exponents $\bm{p}$ and $\bm{q}$ and vectors of coefficients $\bm{\theta}$ and $\bm{\sigma}$ need to be evaluated. We propose to combine alternating computations according to Algorithm~\ref{A-AxeqBy} intended to handle two-sided vector equations, with the approximation technique provided by Algorithm~\ref{A-AP} to solve parameterized one-sided equations that appear in the alternating computations.

We start with the observation that the matrix $\bm{Y}$ is strictly diagonal (does not have zero diagonal entries). We define a matrix $\bm{Y}^{-1}$ obtained from $\bm{Y}$ by replacing all diagonal entries by their inverses, and see that the following equations are equivalent:
\begin{equation*}
\bm{X}(\bm{p})\bm{\theta}
=
\bm{Y}\bm{Z}(\bm{q})\bm{\sigma},
\qquad
\bm{Z}(\bm{q})\bm{\sigma}
=
\bm{Y}^{-1}\bm{X}(\bm{p})\bm{\theta}.
\end{equation*}

\newpage
The alternating computations from Algorithm~\ref{A-AxeqBy} in this context lead to an iterative procedure that alternately solves two equations of the form
\begin{equation*}
\bm{X}(\bm{p})\bm{\theta}
=
\bm{b},
\qquad
\bm{Z}(\bm{q})\bm{\sigma}
=
\bm{a},
\end{equation*}
where $\bm{a}$ and $\bm{b}$ are vectors obtained at previous steps of the procedure. Both equations take the form of \eqref{E-Xpthetaeqy}, and hence we can solve them by using Algorithm~\ref{A-AP}.

As an initial approximation, we set $\bm{\sigma}_{0}=\bm{1}$ and $\bm{q}_{0}=\bm{1}$, and note that then $\bm{Z}(\bm{q}_{0})=\bm{1}\bm{1}^{T}$. Next, we have $\bm{Z}(\bm{q}_{0})\bm{\sigma}_{0}=\bm{1}$, and thus $\bm{Y}\bm{Z}(\bm{q}_{0})\bm{\sigma}_{0}=\bm{y}$.

We denote $\bm{b}_{1}=\bm{Y}\bm{Z}(\bm{q}_{0})\bm{\sigma}_{0}=\bm{y}$ and examine the equation $\bm{X}(\bm{p})\bm{\theta}=\bm{b}_{1}$. This equation is approximated by solving the two-stage minimization problem
\begin{equation*}
\min_{\bm{p}}\min_{\bm{\theta}}\ 
d(\bm{X}(\bm{p})\bm{\theta},\bm{b}_{1}).
\end{equation*}

The inner minimization problem has its minimum given by Theorem~\ref{T-Axeqb} in the form
\begin{equation*}
\delta_{1}(\bm{p})
=
(\bm{X}(\bm{p})(\bm{b}_{1}^{-}\bm{X}(\bm{p}))^{-})^{-}\bm{b}_{1}.
\end{equation*}

Application of Algorithm~\ref{A-AP} to minimize the function $\delta_{1}(\bm{p})$ yields
\begin{equation*}
\bm{p}_{1}
=
\arg\min_{\bm{p}}\ 
\delta_{1}(\bm{p}),
\qquad
\Delta_{1}
=
\delta_{1}(\bm{p}_{1}),
\qquad
\bm{\theta}_{1}
=
\sqrt{\Delta_{1}}
(\bm{b}_{1}^{-}\bm{X}(\bm{p}_{1}))^{-}.
\end{equation*}

Furthermore, we calculate the vector $\bm{a}_{2}=\bm{Y}^{-1}\bm{X}(\bm{p}_{1})\bm{\theta}_{1}$ and then approximate the equation $\bm{Z}(\bm{q})\bm{\sigma}=\bm{a}_{2}$. As before, we represent the approximation problem as follows:
\begin{equation*}
\min_{\bm{q}}\min_{\bm{\sigma}}\ 
d(\bm{Z}(\bm{q})\bm{\sigma},\bm{a}_{2}).
\end{equation*}

According to Theorem~\ref{T-Axeqb}, the minimum of the inner problem is written as
\begin{equation*}
\delta_{2}(\bm{q})
=
(\bm{Z}(\bm{q})(\bm{a}_{2}^{-}\bm{Z}(\bm{q}))^{-})^{-}\bm{a}_{2}.
\end{equation*}

We use Algorithm~\ref{A-AP} to minimize $\delta_{2}(\bm{q})$, which gives
\begin{equation*}
\bm{q}_{2}
=
\arg\min_{\bm{q}}\delta_{2}(\bm{q}),
\qquad
\Delta_{2}
=
\delta_{2}(\bm{q}_{2}),
\qquad
\bm{\sigma}_{2}
=
\sqrt{\Delta_{2}}
(\bm{a}_{2}^{-}\bm{Z}(\bm{q}_{2}))^{-}.
\end{equation*}

In the next two iterations, we first calculate the vector $\bm{b}_{3}=\bm{Y}\bm{Z}(\bm{q}_{2})\bm{\sigma}_{2}$ and approximate the equation $\bm{X}(\bm{p})\bm{\theta}=\bm{b}_{3}$ to find the vectors $\bm{p}_{3}$ and $\bm{\theta}_{3}$. Then, we calculate $\bm{a}_{4}=\bm{Y}^{-1}\bm{X}(\bm{p}_{3})\bm{\theta}_{3}$ and approximate the equation $\bm{Z}(\bm{q})\bm{\sigma}=\bm{a}_{4}$ to find $\bm{q}_{4}$ and $\bm{\sigma}_{4}$.

The iterations continue while the squared error $\Delta_{k}$ is decreasing and stop when the change in two successive squared errors is within a predefined tolerance.  

We give a formal description of the procedure in the form of Algorithm~\ref{A-ARF}.

\vspace{+12pt}

\begin{algorithm}[H]
\caption{{Approximation}
 by rational functions in max-plus algebra.}
\label{A-ARF}
\begin{enumerate}
\item
Input $M$ samples $(x_{i},y_{i})$ for $i=1,\ldots,M$, and fix positive integers $N$ and $L$.

Define the matrices and vectors
\begin{gather*}
\bm{X}(\bm{p})
=
\left(
\begin{array}{ccc}
x_{1}^{p_{1}} & \ldots & x_{1}^{p_{N}}
\\
\vdots & & \vdots
\\
x_{M}^{p_{1}} & \ldots & x_{M}^{p_{N}}
\end{array}
\right),
\qquad
\bm{Y}
=
\left(
\begin{array}{ccc}
y_{1} &  & \mathbb{0}
\\
& \ddots &
\\
\mathbb{0} & & y_{M}
\end{array}
\right),
\\
\bm{Z}(\bm{q})
=
\left(
\begin{array}{ccc}
x_{1}^{q_{1}} & \ldots & x_{1}^{q_{L}}
\\
\vdots & & \vdots
\\
x_{M}^{q_{1}} & \ldots & x_{M}^{q_{L}}
\end{array}
\right),
\qquad
\bm{q}_{0}
=
\bm{1},
\qquad
\bm{\sigma}_{0}
=
\bm{1}.
\end{gather*}
Set a squared error tolerance $\varepsilon>0$, and put $k=1$.

\item
\label{S-Calculate_bk}
Calculate the vector and define the function
\begin{equation*}
\bm{b}_{k}
=
\bm{Y}
\bm{Z}(\bm{q}_{k-1})\bm{\sigma}_{k-1},
\qquad
\delta_{k}(\bm{p})
=
(\bm{X}(\bm{p})(\bm{b}_{k}^{-}\bm{X}(\bm{p}))^{-})^{-}\bm{b}_{k}.
\end{equation*}
\item
Apply Algorithm~\ref{A-AP} to minimize the function $\delta_{k}(\bm{p})$, which yields
\begin{equation*}
\bm{p}_{k}
=
\arg\min_{\bm{p}}
\delta_{k}(\bm{p}),
\qquad
\Delta_{k}
=
\delta_{k}(\bm{p}_{k}),
\qquad
\bm{\theta}_{k}
=
\sqrt{\Delta_{k}}
(\bm{b}_{k}^{-}\bm{X}(\bm{p}_{k}))^{-}.
\end{equation*}
\item
Set $k$ to $k+1$, and define the vector and function
\begin{equation*}
\bm{a}_{k}
=
\bm{Y}^{-1}
\bm{X}(\bm{p}_{k-1})\bm{\theta}_{k-1},
\qquad
\delta_{k}(\bm{q})
=
(\bm{Z}(\bm{q})(\bm{a}_{k}^{-}\bm{Z}(\bm{q}))^{-})^{-}\bm{a}_{k}.
\end{equation*}
\item
Apply Algorithm~\ref{A-AP} to minimize the function $\delta_{k}(\bm{q})$, which yields
\begin{equation*}
\bm{q}_{k}
=
\arg\min_{\bm{q}}
\delta_{k}(\bm{q}),
\qquad
\Delta_{k}
=
\delta_{k}(\bm{q}_{k}),
\qquad
\bm{\sigma}_{k}
=
\sqrt{\Delta_{k}}
(\bm{a}_{k}^{-}\bm{Z}(\bm{q}_{k}))^{-}.
\end{equation*}
\item
If $\Delta_{k}>\Delta_{k-1}$, then put
\begin{equation*}
\Delta_{\ast}
=
\Delta_{k-1},
\qquad
\bm{p}_{\ast}
=
\bm{p}_{k-1},
\qquad
\bm{\theta}_{\ast}
=
\bm{\theta}_{k-1},
\qquad
\bm{q}_{\ast}
=
\bm{q}_{k-2},
\qquad
\bm{\sigma}_{\ast}
=
\bm{\sigma}_{k-2},
\end{equation*}
and stop; otherwise, if $\Delta_{k}>\varepsilon^{-1}\Delta_{k-1}$, then put
\begin{equation*}
\Delta_{\ast}
=
\Delta_{k},
\qquad
\bm{p}_{\ast}
=
\bm{p}_{k-1},
\qquad
\bm{\theta}_{\ast}
=
\bm{\theta}_{k-1},
\qquad
\bm{q}_{\ast}
=
\bm{q}_{k},
\qquad
\bm{\sigma}_{\ast}
=
\bm{\sigma}_{k},
\end{equation*}
and stop; otherwise set $k$ to $k+1$, and go to Step~\ref{S-Calculate_bk}.
\end{enumerate}
\end{algorithm}
\vspace{+12pt}

In the same way as for Algorithm~\ref{A-AP}, the above procedure can be readily adapted to solve rational approximation problems in terms of max-algebra.

\section{Numerical Examples and Graphical Illustrations}
\label{S-NEGI}

In this section, we offer examples to present numerical results of the solution of a best discrete approximation problem where a real nonconvex function is approximated in the max-plus algebra setting. We consider a problem of fitting max-plus rational functions to sample data obtained from a function defined in terms of conventional algebra as follows: 
\begin{equation*}
f(x)
=
3(x-1)^{2}\sin(x)+1/4,
\qquad
x\in[0,2].
\end{equation*}

Given $M=21$ values of the input $x_{i}=1+(i-1)/10$ and output $y_{i}=f(x_{i})$ of the function for $i=1,\ldots,M$, the problem is to approximate $f(x)$ by a tropical rational function
\begin{equation*}
R(x)
=
P(x)/Q(x),
\end{equation*}
defined as a ratio of max-plus algebra polynomials with $N$ and $L$ monomials in the form
\begin{equation*}
P(x)
=
\theta_{1}x^{p_{1}}\oplus\cdots\oplus\theta_{N}x^{p_{N}},
\qquad
Q(x)
=
\sigma_{1}x^{q_{1}}\oplus\cdots\oplus\sigma_{L}x^{q_{L}}.
\end{equation*}

The solution involves the evaluation of vectors $\bm{p}_{\ast}=(p_{1}^{\ast},\ldots,p_{N}^{\ast})^{T}$, $\bm{\theta}_{\ast}=(\theta_{1}^{\ast},\ldots,\theta_{N}^{\ast})^{T}$, $\bm{q}_{\ast}=(q_{1}^{\ast},\ldots,q_{L}^{\ast})^{T}$ and $\bm{\sigma}_{\ast}=(\sigma_{1}^{\ast},\ldots,\sigma_{L}^{\ast})^{T}$ that determine the approximating polynomials $P_{\ast}(x)$ and $Q_{\ast}(x)$ and hence the approximating function $R_{\ast}(x)=P_{\ast}(x)/Q_{\ast}(x)$.

To handle the problem, we apply Algorithm~\ref{A-ARF} under different settings of the numbers $N$ and $L$ to find the best squared approximation error $\Delta_{\ast}$ together with approximating vectors of exponents and coefficients. The tolerance for the squared error is set to $\varepsilon=0.0001$.

The results obtained include numerical solutions for the problems where the number $N$ of monomials for the polynomial $P(x)$ and the number $L$ for $Q(x)$ take successive values from 2 to $N$. Below, for each fixed $N$, we demonstrate the solution corresponding to that $L$ which yields the minimal approximation error. 

To conduct numerical experiments, we use a MATLAB (Release R2021a) code, which combines a collection of functions for the basic computations of tropical algebra with functions that implement the main functionality of the procedure. The code is run on a custom computer equipped with a four-core eight-thread Intel Xeon E3-1231 v3 CPU at 3.40GHz and 32GB of DDR3 RAM, and operated under Windows~10 Enterprise 64-bit OS.  

We start by setting $N=2$, which leads to a minimal squared approximation error $\Delta_{\ast}=0.3099$ that is attained for $L=2$ (the results obtained for $L>2$ do not improve the error). The approximating vectors of exponents and coefficients are given by
\begin{gather*}
\bm{p}_{\ast}
=
\left(
\begin{array}{r}
-0.0628
\\
3.8735
\end{array}
\right),
\quad
\bm{\theta}_{\ast}
=
\left(
\begin{array}{r}
0.5100
\\
-4.6017
\end{array}
\right),
\quad
\bm{q}_{\ast}
=
\left(
\begin{array}{r}
-2.4888
\\
0.1216
\end{array}
\right),
\quad
\bm{\sigma}_{\ast}
=
\left(
\begin{array}{r}
0.4150
\\
-0.0793
\end{array}
\right).
\end{gather*}

The approximating function is written in terms of conventional algebra as
\begin{multline*}
R_{\ast}(x)
=
\max(-0.0628x+0.5100, 3.8735x-4.6017)
\\
-
\max(-2.4888x+0.4150, 0.1216-0.0793).
\end{multline*}

\begin{figure}[H]
\begin{tikzpicture}

\begin{axis}[legend pos=south east,
width=14cm,
height=8cm,
grid=major,
ymin=0,
ymax=3,
xmin=0,
xmax=2
]

\addplot[
only marks,
]
coordinates {
(0.0000,0.2500)
(0.1000,0.4926)
(0.2000,0.6314)
(0.3000,0.6844)
(0.4000,0.6706)
(0.5000,0.6096)
(0.6000,0.5210)
(0.7000,0.4239)
(0.8000,0.3361)
(0.9000,0.2735)
(1.0000,0.2500)
(1.1000,0.2767)
(1.2000,0.3618)
(1.3000,0.5102)
(1.4000,0.7230)
(1.5000,0.9981)
(1.6000,1.3295)
(1.7000,1.7077)
(1.8000,2.1198)
(1.9000,2.5495)
(2.0000,2.9779)
};

\addlegendentry{$y_{i}$}

\addplot[
blue,
samples=200,
line width=1.75pt,
domain=0.0:2.0,
y domain=0:3
]{max(
 -0.0628*x+0.5100,
  3.8735*x-4.6017   
)
-max(
 -2.4888*x+0.4150,
  0.1216*x-0.0793
)
};

\addlegendentry{$R_{\ast}(x)$}

\addplot[
samples=200,
black,
line width=0.75pt,
domain=0.0:2.0,
y domain=0:3
]{3*(x-1)^2*sin(deg(x))+1/4};

\addlegendentry{$f(x)$}

\node[style={fill=white}] at (axis cs: 1.0,2.5) {$f(x)=3(x-1)^{2}\sin(x)+1/4$};
\node[style={fill=white}] at (axis cs: 1.0,2.0) {$\Delta_{\ast}=0.3099$};

\end{axis}

\end{tikzpicture}
\caption{Approximation by a max-plus rational function $R_{\ast}(x)$ with $N=2$ and $L=2$.}
\label{F-AMPRFN2L2}
\end{figure}
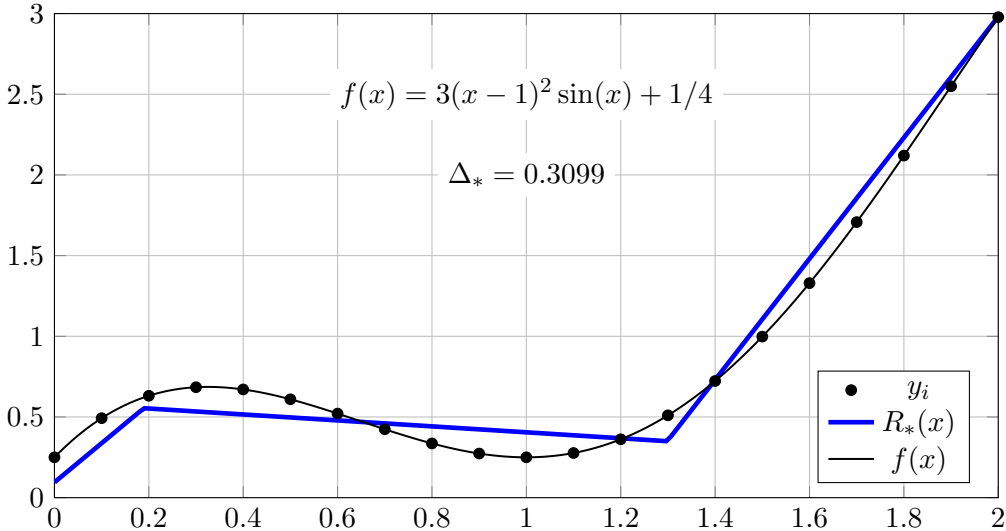

A graphical representation of the solution is shown in Figure~\ref{F-AMPRFN2L2}.


Furthermore, we set the number $N$ of monomials in the polynomial $P(x)$ to $3$ and solve the approximation problem for various numbers $L$ of monomials in $Q(x)$. The minimal squared approximation error is achieved if $L=3$ to be $\Delta_{\ast}=0.1158$. The approximation procedure yields vectors of exponents and coefficients that take the form
\begin{equation*}
\bm{p}_{\ast}
=
\left(
\begin{array}{r}
-0.5398
\\
2.2198
\\
4.3520
\end{array}
\right),
\hspace*{2mm}
\bm{\theta}_{\ast}
=
\left(
\begin{array}{r}
0.8577
\\
-2.2300
\\
-5.5099
\end{array}
\right),
\quad
\bm{q}_{\ast}
=
\left(
\begin{array}{r}
-1.9879
\\
0.0228
\\
0.2545
\end{array}
\right),
\quad
\bm{\sigma}_{\ast}
=
\left(
\begin{array}{r}
0.5678
\\
0.0200
\\
-0.2349
\end{array}
\right).
\end{equation*}

\begin{figure}[H]
\begin{tikzpicture}

\begin{axis}[legend pos=south east,
width=14cm,
height=8cm,
grid=major,
ymin=0,
ymax=3,
xmin=0,
xmax=2
]

\addplot[
only marks,
]
coordinates {
(0.0000,0.2500)
(0.1000,0.4926)
(0.2000,0.6314)
(0.3000,0.6844)
(0.4000,0.6706)
(0.5000,0.6096)
(0.6000,0.5210)
(0.7000,0.4239)
(0.8000,0.3361)
(0.9000,0.2735)
(1.0000,0.2500)
(1.1000,0.2767)
(1.2000,0.3618)
(1.3000,0.5102)
(1.4000,0.7230)
(1.5000,0.9981)
(1.6000,1.3295)
(1.7000,1.7077)
(1.8000,2.1198)
(1.9000,2.5495)
(2.0000,2.9779)
};

\addlegendentry{$y_{i}$}

\addplot[
blue,
samples=200,
line width=1.75pt,
domain=0.0:2.0,
y domain=0:3
]{max(
 -0.5398*x+0.8577,
  2.2198*x-2.2300,
  4.3520*x-5.5099
)
-max(
 -1.9879*x+0.5678,
  0.0228*x+0.0200,
  0.2545*x-0.2349
)
};

\addlegendentry{$R_{\ast}(x)$}

\addplot[
samples=200,
black,
line width=0.75pt,
domain=0.0:2.0,
y domain=0:3
]{3*(x-1)^2*sin(deg(x))+1/4};

\addlegendentry{$f(x)$}

\node[style={fill=white}] at (axis cs: 1.0,2.5) {$f(x)=3(x-1)^{2}\sin(x)+1/4$};
\node[style={fill=white}] at (axis cs: 1.0,2.0) {$\Delta_{\ast}=0.1158$};

\end{axis}

\end{tikzpicture}
\caption{Approximation by a max-plus rational function $R_{\ast}(x)$ with $N=3$ and $L=3$.}
\label{F-AMPRFN3L3}
\end{figure}
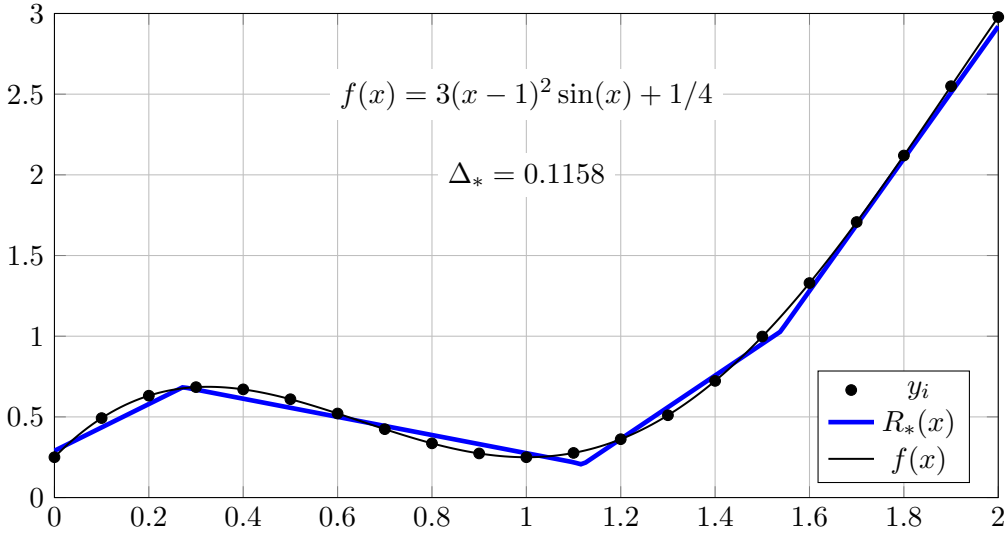

\newpage
The obtained approximating function is given using standard arithmetic operations as
\begin{multline*}
R_{\ast}(x)
=
\max(-0.5398x+0.8577, 2.2198x-2.2300, 4.3520x-5.5099)
\\
-
\max(-1.9879x+0.5678, 0.0228x+0.0200, 0.2545x-0.2349),
\end{multline*}
and graphically illustrated in Figure~\ref{F-AMPRFN3L3}.

\vspace{-6pt}

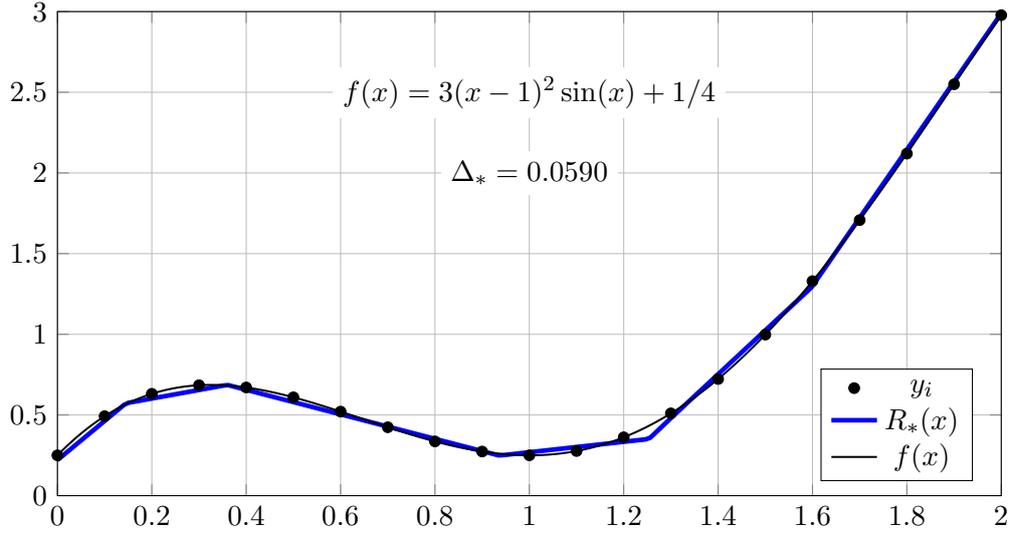
\begin{figure}[H]
\begin{tikzpicture}

\begin{axis}[legend pos=south east,
width=14cm,
height=8cm,
grid=major,
ymin=0,
ymax=3,
xmin=0,
xmax=2
]

\addplot[
only marks,
]
coordinates {
(0.0000,0.2500)
(0.1000,0.4926)
(0.2000,0.6314)
(0.3000,0.6844)
(0.4000,0.6706)
(0.5000,0.6096)
(0.6000,0.5210)
(0.7000,0.4239)
(0.8000,0.3361)
(0.9000,0.2735)
(1.0000,0.2500)
(1.1000,0.2767)
(1.2000,0.3618)
(1.3000,0.5102)
(1.4000,0.7230)
(1.5000,0.9981)
(1.6000,1.3295)
(1.7000,1.7077)
(1.8000,2.1198)
(1.9000,2.5495)
(2.0000,2.9779)
};

\addlegendentry{$y_{i}$}

\addplot[
blue,
samples=200,
line width=1.75pt,
domain=0.0:2.0,
y domain=0:3
]{max(
 -0.6204*x+0.8898,
  0.4634*x-0.1250,
  2.9159*x-3.1989,
  4.4184*x-5.5992
)
-max(
 -3.0464*x+0.6693,
 -1.1501*x+0.3938,
  0.1395*x-0.0702,
  0.1971*x-0.1452
)
};

\addlegendentry{$R_{\ast}(x)$}

\addplot[
samples=200,
black,
line width=0.75pt,
domain=0.0:2.0,
y domain=0:3
]{3*(x-1)^2*sin(deg(x))+1/4};

\addlegendentry{$f(x)$}

\node[style={fill=white}] at (axis cs: 1.0,2.5) {$f(x)=3(x-1)^{2}\sin(x)+1/4$};
\node[style={fill=white}] at (axis cs: 1.0,2.0) {$\Delta_{\ast}=0.0590$};

\end{axis}

\end{tikzpicture}
\caption{Approximation by a max-plus rational function $R_{\ast}(x)$ with $N=4$ and $L=4$.}
\label{F-AMPRFN4L4}
\end{figure}

In the case when $N=4$, the minimal squared error is found for $L=4$, and it is equal to $\Delta_{\ast}=0.0590$. The obtained vectors of exponents and coefficients are 
\begin{equation*}
\bm{p}_{\ast}
=
\left(
\begin{array}{r}
-0.6204
\\
0.4634
\\
2.9159
\\
4.4184
\end{array}
\right),
\hspace*{2mm}
\bm{\theta}_{\ast}
=
\left(
\begin{array}{r}
0.8898
\\
-0.1250
\\
-3.1989
\\
-5.5992
\end{array}
\right),
\quad
\bm{q}_{\ast}
=
\left(
\begin{array}{r}
-3.0464
\\
-1.1501
\\
0.1395
\\
0.1971
\end{array}
\right),
\quad
\bm{\sigma}_{\ast}
=
\left(
\begin{array}{r}
0.6693
\\
0.3938
\\
-0.0702
\\
-0.1452 
\end{array}
\right),
\end{equation*}
which determines the approximating function (see Figure~\ref{F-AMPRFN4L4}) in the form
\begin{multline*}
R_{\ast}(x)
=
\max(-0.6204x+0.8898, 0.4634x-0.1250, 2.9159x-3.1989, 4.4184x-5.5992)
\\
-
\max(-3.0464x+0.6693, -1.1501x+0.3938, 0.1395x-0.0702, 0.1971x-0.1452).
\end{multline*}\vspace{-25pt}

Under the condition that $N=5$, the minimal squared error $\Delta_{\ast}=0.0370$ is attained when $L=3$. The obtained vectors of exponents and coefficients are given by
\begin{equation*}
\bm{p}_{\ast}
=
\left(
\begin{array}{r}
-0.6350
\\
0.5005
\\
1.7163
\\
2.9842
\\
4.3541
\end{array}
\right),
\hspace*{2mm}
\bm{\theta}_{\ast}
=
\left(
\begin{array}{r}
0.8840
\\
-0.1512
\\
-1.5518
\\
-3.2623
\\
-5.4317 
\end{array}
\right),
\quad
\bm{q}_{\ast}
=
\left(
\begin{array}{r}
-3.0609
\\
-1.1647
\\
0.2332 
\end{array}
\right),
\quad
\bm{\sigma}_{\ast}
=
\left(
\begin{array}{r}
0.6525
\\
0.3770
\\
-0.1524 
\end{array}
\right).
\end{equation*}

The approximating function takes the form
\begin{multline*}
R_{\ast}(x)
=
\max(-0.6350x+0.8840, 0.5005x-0.1512, 1.7163x-1.5518,
\\
2.9842x-3.2623, 4.3541x-5.4317)
\\
-\max(-3.0609x+0.6525, -1.1647x+0.3770, 0.2332x-0.1524).
\end{multline*}

The solution obtained is shown in Figure~\ref{F-AMPRFN5L3}.


In the case when $N=6$, the minimal squared error is $\Delta_{\ast}=0.0113$, which is achieved for $L=5$. The approximating vectors have the form
\begin{equation*}
\bm{p}_{\ast}
=
\left(
\begin{array}{r}
-0.6251
\\
0.0174
\\
1.1035
\\
2.3810
\\
3.5666
\\
4.4863
\end{array}
\right),
\hspace*{2mm}
\bm{\theta}_{\ast}
=
\left(
\begin{array}{r}
0.8792
\\
0.3250
\\
-0.8195
\\
-2.4170
\\
-4.1332
\\
-5.6611 
\end{array}
\right),
\quad
\bm{q}_{\ast}
=
\left(
\begin{array}{r}
-3.0511
\\
-1.1548
\\
-0.0151
\\
0.2506
\\
0.2524
\end{array}
\right),
\quad
\bm{\sigma}_{\ast}
=
\left(
\begin{array}{r}
0.6348
\\
0.3593
\\
-0.0298
\\
-0.1626
\\
-0.1657 
\end{array}
\right).
\end{equation*}

\begin{figure}[H]
\begin{tikzpicture}

\begin{axis}[legend pos=south east,
width=14cm,
height=8cm,
grid=major,
ymin=0,
ymax=3,
xmin=0,
xmax=2
]

\addplot[
only marks,
]
coordinates {
(0.0000,0.2500)
(0.1000,0.4926)
(0.2000,0.6314)
(0.3000,0.6844)
(0.4000,0.6706)
(0.5000,0.6096)
(0.6000,0.5210)
(0.7000,0.4239)
(0.8000,0.3361)
(0.9000,0.2735)
(1.0000,0.2500)
(1.1000,0.2767)
(1.2000,0.3618)
(1.3000,0.5102)
(1.4000,0.7230)
(1.5000,0.9981)
(1.6000,1.3295)
(1.7000,1.7077)
(1.8000,2.1198)
(1.9000,2.5495)
(2.0000,2.9779)
};

\addlegendentry{$y_{i}$}

\addplot[
blue,
samples=200,
line width=1.75pt,
domain=0.0:2.0,
y domain=0:3
]{max(
 -0.6251*x+0.8792,
  0.0174*x+0.3250,
  1.1035*x-0.8195,
  2.3810*x-2.4170,
  3.5666*x-4.1332,
  4.4863*x-5.6611
)
-max(
 -3.0511*x+0.6348,
 -1.1548*x+0.3593,
 -0.0151*x-0.0298,
  0.2506*x-0.1626,
  0.2524*x-0.1657
)
};

\addlegendentry{$R_{\ast}(x)$}

\addplot[
samples=200,
black,
line width=0.75pt,
domain=0.0:2.0,
y domain=0:3
]{3*(x-1)^2*sin(deg(x))+1/4};

\addlegendentry{$f(x)$}

\node[style={fill=white}] at (axis cs: 1.0,2.5) {$f(x)=3(x-1)^{2}\sin(x)+1/4$};
\node[style={fill=white}] at (axis cs: 1.0,2.0) {$\Delta_{\ast}=0.0113$};

\end{axis}

\end{tikzpicture}
\caption{Approximation by a max-plus rational function $R_{\ast}(x)$ with $N=6$ and $L=5$.}
\label{F-AMPRFN6L5}
\end{figure}
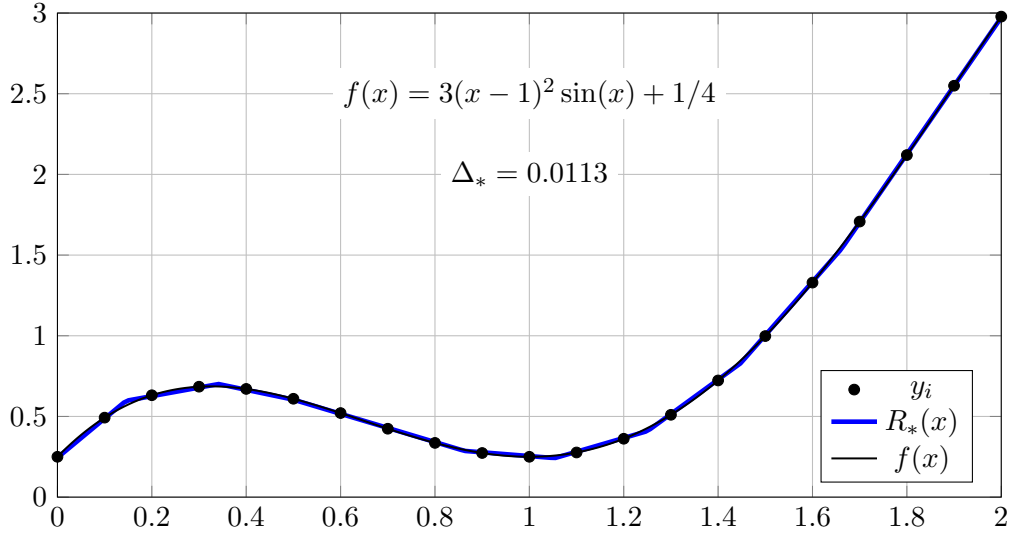

The approximating function (see Figure~\ref{F-AMPRFN6L5}) is written as
\begin{multline*}
R_{\ast}(x)
=
\max(-0.6251x+0.8792, 0.0174x+0.3250, 1.1035x-0.8195, 2.3810x-2.4170,
\\
3.5666x-4.1332, 4.4863x-5.6611),
\\
-
\max(-3.0511x+0.6348, -1.1548x+0.3593, -0.0151x-0.0298,
\\
0.2506x-0.1626, 0.2524x-0.1657).
\end{multline*}\vspace{-28pt}

\begin{figure}[H]
\begin{tikzpicture}

\begin{axis}[legend pos=south east,
width=14cm,
height=8cm,
grid=major,
ymin=0,
ymax=3,
xmin=0,
xmax=2
]

\addplot[
only marks,
]
coordinates {
(0.0000,0.2500)
(0.1000,0.4926)
(0.2000,0.6314)
(0.3000,0.6844)
(0.4000,0.6706)
(0.5000,0.6096)
(0.6000,0.5210)
(0.7000,0.4239)
(0.8000,0.3361)
(0.9000,0.2735)
(1.0000,0.2500)
(1.1000,0.2767)
(1.2000,0.3618)
(1.3000,0.5102)
(1.4000,0.7230)
(1.5000,0.9981)
(1.6000,1.3295)
(1.7000,1.7077)
(1.8000,2.1198)
(1.9000,2.5495)
(2.0000,2.9779)
};

\addlegendentry{$y_{i}$}

\addplot[
blue,
samples=200,
line width=1.75pt,
domain=0.0:2.0,
y domain=0:3
]{max(
 -0.6350*x+0.8840,
  0.5005*x-0.1512,
  1.7163*x-1.5518,
  2.9842*x-3.2623,
  4.3541*x-5.4317
)
-max(
 -3.0609*x+0.6525,
 -1.1647*x+0.3770,
  0.2332*x-0.1524
)
};

\addlegendentry{$R_{\ast}(x)$}

\addplot[
samples=200,
black,
line width=0.75pt,
domain=0.0:2.0,
y domain=0:3
]{3*(x-1)^2*sin(deg(x))+1/4};

\addlegendentry{$f(x)$}

\node[style={fill=white}] at (axis cs: 1.0,2.5) {$f(x)=3(x-1)^{2}\sin(x)+1/4$};
\node[style={fill=white}] at (axis cs: 1.0,2.0) {$\Delta_{\ast}=0.0370$};

\end{axis}

\end{tikzpicture}
\caption{Approximation by a max-plus rational function $R_{\ast}(x)$ with $N=5$ and $L=3$.}
\label{F-AMPRFN5L3}
\end{figure}
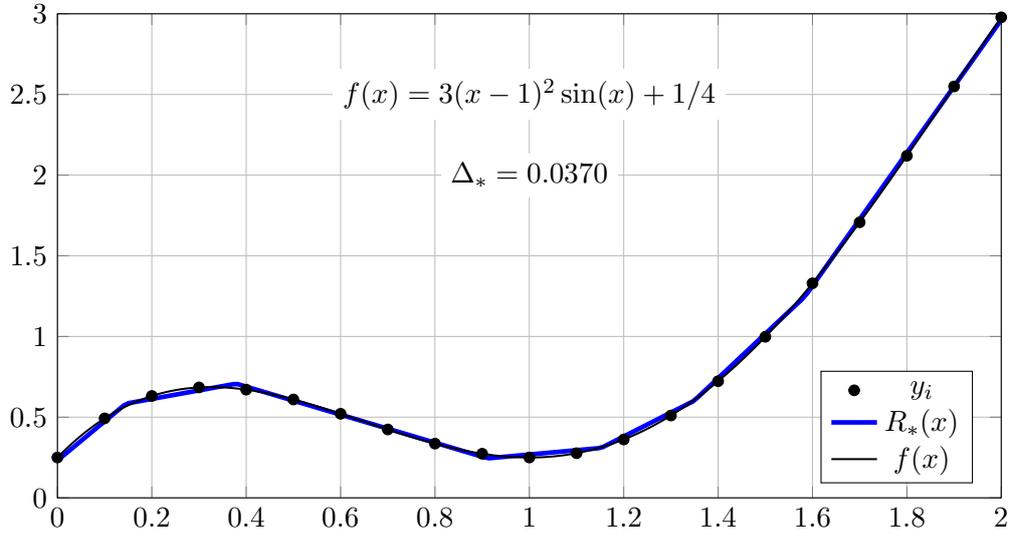

Finally, for $N=7$, we obtain the minimal error $\Delta_{\ast}<\varepsilon$ if $L=5$. We observe that in this case, the solution of the approximation problem leads to the exact solution with $\Delta_{\ast}=0$. The obtained vectors of exponents and coefficients are as follows:
\begin{equation*}
\bm{p}_{\ast}
=
\left(
\begin{array}{r}
-0.6963
\\
-0.3512
\\
0.5419
\\
1.7577
\\
3.0256
\\
4.0567
\\
4.5718
\end{array}
\right),
\hspace*{2mm}
\bm{\theta}_{\ast}
=
\left(
\begin{array}{r}
0.9254
\\
0.6586
\\
-0.1955
\\
-1.5961
\\
-3.3066
\\
-4.9000
\\
-5.7934
\end{array}
\right),
\quad
\bm{q}_{\ast}
=
\left(
\begin{array}{r}
-3.1222
\\
-1.2260
\\
-0.0863
\\
0.2746
\\
0.2880
\end{array}
\right),
\quad
\bm{\sigma}_{\ast}
=
\left(
\begin{array}{r}
0.6754
\\
0.3999
\\
0.0108
\\
-0.1782
\\
-0.2036
\end{array}
\right).
\end{equation*}

The approximating function is given by
\begin{multline*}
R_{\ast}(x)
=
\max(-0.6963x+0.9254, -0.3512x+0.6586, 0.5419x-0.1955, 1.7577x-1.5961,
\\
3.0256x-3.3066, 4.0567x-4.9000, 4.5718x-5.7934),
\\
-
\max(-3.1222x+0.6754, -1.2260x+0.3999, -0.0863x+0.0108,
\\
0.2746x-0.1782, 0.2880x-0.2036).
\end{multline*}

\begin{figure}[H]
\begin{tikzpicture}

\begin{axis}[legend pos=south east,
width=14cm,
height=8cm,
grid=major,
ymin=0,
ymax=3,
xmin=0,
xmax=2
]

\addplot[
only marks,
]
coordinates {
(0.0000,0.2500)
(0.1000,0.4926)
(0.2000,0.6314)
(0.3000,0.6844)
(0.4000,0.6706)
(0.5000,0.6096)
(0.6000,0.5210)
(0.7000,0.4239)
(0.8000,0.3361)
(0.9000,0.2735)
(1.0000,0.2500)
(1.1000,0.2767)
(1.2000,0.3618)
(1.3000,0.5102)
(1.4000,0.7230)
(1.5000,0.9981)
(1.6000,1.3295)
(1.7000,1.7077)
(1.8000,2.1198)
(1.9000,2.5495)
(2.0000,2.9779)
};

\addlegendentry{$y_{i}$}

\addplot[
blue,
samples=200,
line width=1.75pt,
domain=0.0:2.0,
y domain=0:3
]{max(
 -0.6963*x+0.9254,
 -0.3512*x+0.6586,
  0.5419*x-0.1955,
  1.7577*x-1.5961,
  3.0256*x-3.3066,
  4.0567*x-4.9000,
  4.5718*x-5.7934
)
-max(
 -3.1222*x+0.6754,
 -1.2260*x+0.3999,
 -0.0863*x+0.0108,
  0.2746*x-0.1782,
  0.2880*x-0.2036
)
};

\addlegendentry{$R_{\ast}(x)$}

\addplot[
samples=200,
black,
line width=0.75pt,
domain=0.0:2.0,
y domain=0:3
]{3*(x-1)^2*sin(deg(x))+1/4};

\addlegendentry{$f(x)$}

\node[style={fill=white}] at (axis cs: 1.0,2.5) {$f(x)=3(x-1)^{2}\sin(x)+1/4$};
\node[style={fill=white}] at (axis cs: 1.0,2.0) {$\Delta_{\ast}<\varepsilon$};

\end{axis}

\end{tikzpicture}
\caption{Approximation by a max-plus rational function $R_{\ast}(x)$ with $N=7$ and $L=5$.}
\label{F-AMPRFN7L5}
\end{figure}
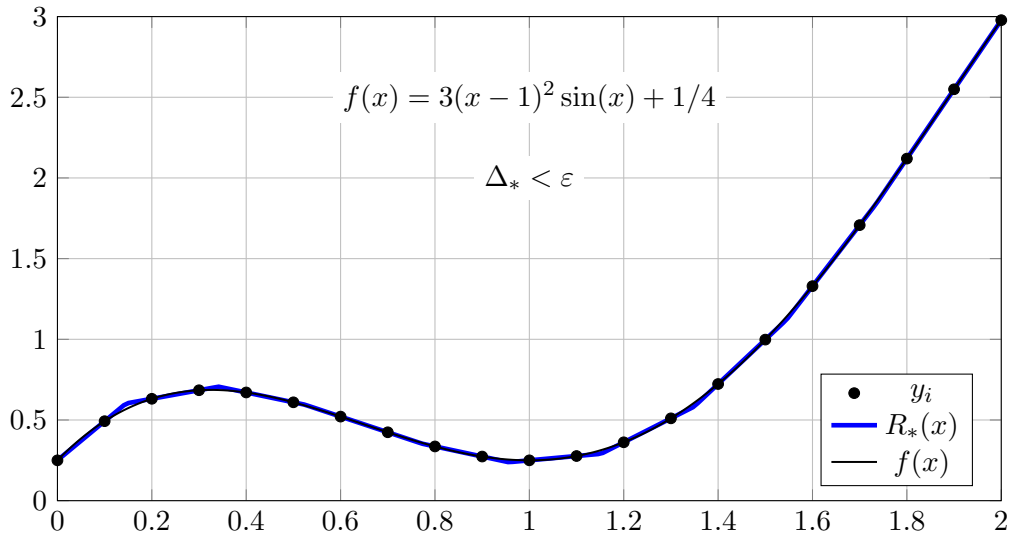

The solution is illustrated in Figure~\ref{F-AMPRFN7L5}.

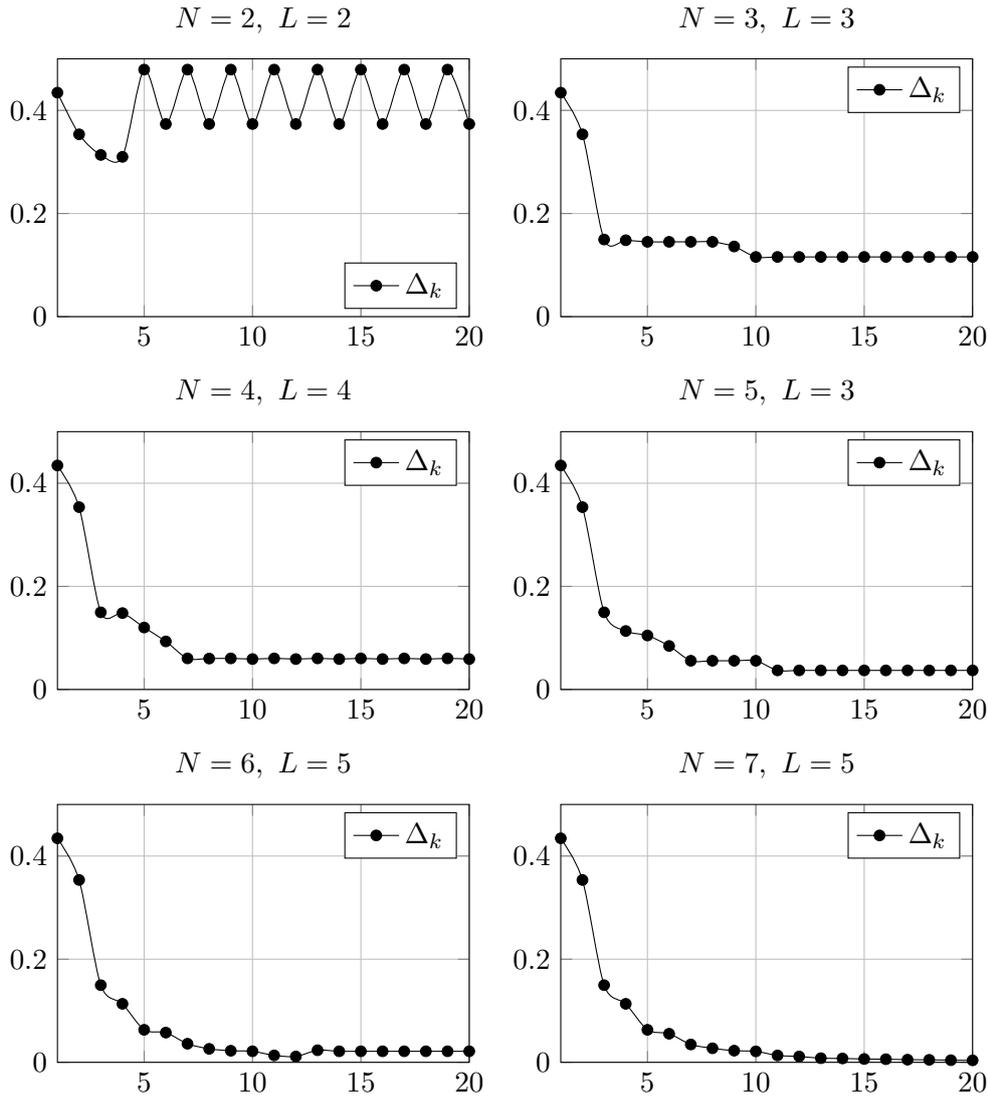
\begin{figure}[H]
\begin{tikzpicture}

\begin{axis}[title={$N=2,\ L=2$},
legend pos=south east,
width=7cm,
height=5cm,
grid=major,
ymin=0,
ymax=0.5,
xmin=1,
xmax=20
]

\addplot[
mark=*,
smooth
]
coordinates {
( 1,0.4344)
( 2,0.3536)
( 3,0.3136)
( 4,0.3099)
( 5,0.4791)
( 6,0.3736)
( 7,0.4791)
( 8,0.3736)
( 9,0.4791)
(10,0.3736)
(11,0.4791)
(12,0.3736)
(13,0.4791)
(14,0.3736)
(15,0.4791)
(16,0.3736)
(17,0.4791)
(18,0.3736)
(19,0.4791)
(20,0.3736)
};

\addlegendentry{$\Delta_{k}$}

\node[style={fill=white}] at (axis cs: 1.7,2.5) {$f(x)=(x-3/4)^{2}-3(x-1)^{1/2}+2$};


\end{axis}

\end{tikzpicture}
\begin{tikzpicture}

\begin{axis}[title={$N=3,\ L=3$},
legend pos=north east,
width=7cm,
height=5cm,
grid=major,
ymin=0,
ymax=0.5,
xmin=1,
xmax=20
]

\addplot[
mark=*,
smooth
]
coordinates {
( 1,0.4344)
( 2,0.3536)
( 3,0.1496)
( 4,0.1482)
( 5,0.1451)
( 6,0.1451)
( 7,0.1451)
( 8,0.1451)
( 9,0.1361)
(10,0.1158)
(11,0.1158)
(12,0.1158)
(13,0.1158)
(14,0.1158)
(15,0.1158)
(16,0.1158)
(17,0.1158)
(18,0.1158)
(19,0.1158)
(20,0.1158)
};

\addlegendentry{$\Delta_{k}$}


\end{axis}

\end{tikzpicture}
\vspace{1ex}

\begin{tikzpicture}

\begin{axis}[title={$N=4,\ L=4$},
legend pos=north east,
width=7cm,
height=5cm,
grid=major,
ymin=0,
ymax=0.5,
xmin=1,
xmax=20
]

\addplot[
mark=*,
smooth
]
coordinates {
( 1,0.4344)
( 2,0.3536)
( 3,0.1496)
( 4,0.1482)
( 5,0.1202)
( 6,0.0932)
( 7,0.0603)
( 8,0.0600)
( 9,0.0603)
(10,0.0590)
(11,0.0603)
(12,0.0590)
(13,0.0603)
(14,0.0590)
(15,0.0603)
(16,0.0590)
(17,0.0603)
(18,0.0590)
(19,0.0603)
(20,0.0590)
};

\addlegendentry{$\Delta_{k}$}


\end{axis}

\end{tikzpicture}
\begin{tikzpicture}

\begin{axis}[title={$N=5,\ L=3$},
legend pos=north east,
width=7cm,
height=5cm,
grid=major,
ymin=0,
ymax=0.5,
xmin=1,
xmax=20
]

\addplot[
mark=*,
smooth
]
coordinates {
( 1,0.4344)
( 2,0.3536)
( 3,0.1496)
( 4,0.1134)
( 5,0.1047)
( 6,0.0844)
( 7,0.0557)
( 8,0.0557)
( 9,0.0557)
(10,0.0557)
(11,0.0370)
(12,0.0370)
(13,0.0370)
(14,0.0370)
(15,0.0370)
(16,0.0370)
(17,0.0370)
(18,0.0370)
(19,0.0370)
(20,0.0370)
};

\addlegendentry{$\Delta_{k}$}


\end{axis}

\end{tikzpicture}
\vspace{1ex}

\begin{tikzpicture}

\begin{axis}[title={$N=6,\ L=5$},
legend pos=north east,
width=7cm,
height=5cm,
grid=major,
ymin=0,
ymax=0.5,
xmin=1,
xmax=20
]

\addplot[
mark=*,
smooth
]
coordinates {
( 1,0.4344)
( 2,0.3536)
( 3,0.1496)
( 4,0.1134)
( 5,0.0630)
( 6,0.0576)
( 7,0.0360)
( 8,0.0261)
( 9,0.0224)
(10,0.0214)
(11,0.0134)
(12,0.0113)
(13,0.0234)
(14,0.0214)
(15,0.0214)
(16,0.0214)
(17,0.0214)
(18,0.0214)
(19,0.0214)
(20,0.0214)
};

\addlegendentry{$\Delta_{k}$}


\end{axis}

\end{tikzpicture}
\begin{tikzpicture}

\begin{axis}[title={$N=7,\ L=5$},
legend pos=north east,
width=7cm,
height=5cm,
grid=major,
ymin=0,
ymax=0.5,
xmin=1,
xmax=20
]

\addplot[
mark=*,
smooth
]
coordinates {
( 1,0.4344)
( 2,0.3536)
( 3,0.1496)
( 4,0.1134)
( 5,0.0630)
( 6,0.0554)
( 7,0.0346)
( 8,0.0272)
( 9,0.0226)
(10,0.0212)
(11,0.0133)
(12,0.0113)
(13,0.0080)
(14,0.0075)
(15,0.0064)
(16,0.0060)
(17,0.0051)
(18,0.0048)
(19,0.0041)
(20,0.0038)
};

\addlegendentry{$\Delta_{k}$}


\end{axis}

\end{tikzpicture}

\caption{Evolution of squared approximation error $\Delta_{k}$ as the number of iterations $k$ increases.}
\label{F-DSAEDI}
\end{figure}

We conclude this section with a graphical representation of the dynamics of the squared approximation error during the iterations. Figure~\ref{F-DSAEDI} demonstrates how the squared error $\Delta_{k}$ changes as the number of iterations $k$ increases in the solutions presented above. 

The graphs in Figure~\ref{F-DSAEDI} show that in all cases except for the first one, the sequence of squared errors $\Delta_{k}$ has a tendency to decrease as $k$ increases. The unusual behavior of the error value in the case when $N=2$ and $L=2$ can be explained by the insufficient number of monomials in the polynomials, which results in unstable behavior of the optimization procedure in Algorithm~\ref{A-MFdeltap} and thus of the overall solution provided by Algorithm~\ref{A-ARF}.

\section{Conclusions}
\label{S-C}

We have considered discrete best approximation problems in the setting of a tropical algebra that deals with semifields and semirings with idempotent addition. Given sample points of the input and output of an unknown function defined on a semifield, the problem is to approximate the function by a rational function in the sense of a generalized metric. The approximating function takes the form of a ratio of Puiseux polynomials that can have real exponents. {The approximation problem requires evaluating both the unknown exponents and coefficients of the monomials in the two approximating polynomials.}

We developed an iterative computational procedure that is based on the transformation into a problem of approximating a vector equation. Each side of the equation corresponds to one of the polynomials and is represented as a matrix parameterized by unknown exponents, multiplied by an unknown vector of coefficients of monomials in the polynomial. The procedure alternately fixes a vector on one side of this equation and solves the new equation for the unknown exponents and coefficients on the other side. The solution of the last equation involves two stages: the first is the approximation of the equation to derive the approximation error and unknown vector parameterized by unknown exponents, and the second is the minimization of the error to find the exponents. To complete the first stage, we exploited the results of approximating a vector equation with a fixed matrix, which are known in analytical form. To handle the minimization problem in the second stage, we applied a combinatorial solution developed in \cite{Krivulin2024Solution}, which is based on agglomerative clustering.

We demonstrated the practical feasibility and computational efficiency of the approximation procedure by numerical and graphical illustrations. In contrast to the algorithm proposed in \cite{Dunbar2024Alternating}, which assumes the exponents in the polynomials to be fixed in advance and concentrates on the evaluation of coefficients, the proposed solution is designed to find the exponents as well as coefficients, which are  both considered unknown.   

We believe the contribution of this study to be twofold. First, the results obtained provide further development of the theory and methods of tropical algebra as an analytical tool to solve tropical approximation problems. The study offers an alternative approach based on new arithmetic to approximation problems, which differs from the classical solutions known in the context of conventional mathematics. Second, the proposed approximation procedure significantly expands the range of applications of tropical algebra in actual real-world problems. The area of applications of the results includes pattern recognition, signal denoising, data compression, neural network design, and other research and applied domains where the best approximation serves as a key instrument of solution.

A possible drawback of the proposed solution is its computational complexity, which may be somewhat higher than that of existing approximation techniques. However, a correct comparison involves a careful formal assessment of the computational complexity of the solution, which requires a separate study. In any case, our numerical experience has shown that the proposed iterative procedure generally converges in a finite number of steps and provides a reasonably good solution to the approximation problems considered. 

Taking into account that the main steps of the procedure have polynomial complexity, and considering the experience gained from the numerical solution of test problems, it can be expected that the proposed procedure is scalable. We reckon that the procedure is quite capable of solving large-scale problems in a reasonable time, but more substantiated conclusions require further experiments, which are intended for future research.

The combinatorial algorithm based on agglomerative clustering, proposed for solving optimization problems, presents the most time-consuming component of the solution. Although this algorithm demonstrates sufficient efficiency in practice, the question of using other approaches to speed up the solution (heuristic algorithms, hybrid methods, parallel computation, etc.) remains open and requires additional investigation.

As was shown in \cite{Krivulin2024Solution}, the approximation procedure developed in the framework of max-plus algebra can be readily adjusted to handle approximation problems in terms of max-algebra. Since the rational functions in the context of max-algebra take the form of splines, the obtained results can be applied for spline approximation of real functions, which needs further experimental and theoretical analysis of the procedure.

The numerical examples given as illustrations were concentrated on a smooth nonconvex function. However, we expect that the proposed procedure will prove to be quite efficient in problems of approximating other types of functions. To demonstrate the robustness and versatility of the approach, further research is necessary to consider approximation problems for a wide range of functions, including non-smooth and discontinuous functions.

Other areas of future research include a formal analysis of the convergence of the procedure and the estimation of its computational complexity. A comprehensive analytical and numerical comparison with existing solutions presents another urgent problem to address.

\bibliographystyle{abbrvurl}

\bibliography{Tropical_solution_of_discrete_best_approximation_problems}

\begin{thebibliography}{10}

\bibitem{Akian2011Best}
M.~Akian, S.~Gaubert, V.~Ni\c{t}ic\u{a}, and I.~Singer.
\newblock Best approximation in max-plus semimodules.
\newblock {\em Linear Algebra Appl.}, 435(12):3261--3296, 2011.
\newblock \href {https://doi.org/10.1016/j.laa.2011.06.009}
  {\path{doi:10.1016/j.laa.2011.06.009}}.

\bibitem{Bereau2024Tropical}
A.~B\'{e}reau.
\newblock {\em Tropical polynomial systems and game theory}.
\newblock PhD thesis, Institut Polytechnique de Paris, 2024.

\bibitem{Butkovic2010Maxlinear}
P.~Butkovi\v{c}.
\newblock {\em Max-linear Systems}.
\newblock Springer Monographs in Mathematics. Springer, London, 2010.
\newblock \href {https://doi.org/10.1007/978-1-84996-299-5}
  {\path{doi:10.1007/978-1-84996-299-5}}.

\bibitem{Cameron1966Piecewise}
S.~H. Cameron.
\newblock Piece-wise linear approximation.
\newblock Technical Note CSTN-106, Computer Sciences Division, IIT Research
  Institute, Chicago, IL, 1966.

\bibitem{Camponogara2015Models}
E.~Camponogara and L.~F. Nazari.
\newblock Models and algorithms for optimal piecewise-linear function
  approximation.
\newblock {\em Math. Probl. Eng.}, 2015:876862, 2015.
\newblock \href {https://doi.org/10.1155/2015/876862}
  {\path{doi:10.1155/2015/876862}}.

\bibitem{Chen2024Tropical}
J.~Chen, D.~Grigoriev, and V.~Shpilrain.
\newblock Tropical cryptography {III}: digital signature.
\newblock {\em J. Math. Cryptol.}, 18(1):20240005, 2024.
\newblock \href {https://doi.org/10.1515/jmc-2024-0005}
  {\path{doi:10.1515/jmc-2024-0005}}.

\bibitem{Conn1988Computational}
A.~R. Conn and Y.~Li.
\newblock The computational structure and characterization of nonlinear
  discrete {C}hebyshev problem.
\newblock Technical Report 88-956, Department of Computer Science, Cornell
  University, Ithaca, NY, 1988.

\bibitem{Laplace1832Mecanique}
P.~S. de~Laplace.
\newblock {\em M{\'e}canique C{\'e}leste. Volume 2}.
\newblock Hillard, Gray, Littl{\`e}, and Wilkins, Boston, 1832.
\newblock (Engl. transl. with comment. by N.~Bowditch).

\bibitem{Dellaccio2025Truncated}
F.~Dell{'A}ccio, A.~Guessab, G.~V. Milovanovi\'{c}, and F.~Nudo.
\newblock Truncated {G}egenbauer-{H}ermite weighted approach for the enrichment
  of the {C}rouzeix-{R}aviart finite element.
\newblock {\em BIT Numer. Math.}, 65(2):24, 2025.
\newblock \href {https://doi.org/10.1007/s10543-025-01069-6}
  {\path{doi:10.1007/s10543-025-01069-6}}.

\bibitem{Dunbar2024Alternating}
A.~Dunbar and L.~Ruthotto.
\newblock Alternating minimization for regression with tropical rational
  functions.
\newblock {\em Alg. Stat.}, 15(1):85--111, 2024.
\newblock \href {https://doi.org/10.2140/astat.2024.15.85}
  {\path{doi:10.2140/astat.2024.15.85}}.

\bibitem{Durcheva2025Closer}
M.~Durcheva.
\newblock Tropical cryptography -- the state of the art and future prospects.
\newblock {\em Athens J. Sci.}, 12(3):189--204, 2025.
\newblock \href {https://doi.org/10.30958/ajs.12-3-3}
  {\path{doi:10.30958/ajs.12-3-3}}.

\bibitem{Esparza2008Approximative}
J.~Esparza, T.~Gawlitza, S.~Kiefer, and H.~Seidl.
\newblock Approximative methods for monotone systems of min-max-polynomial
  equations.
\newblock In L.~Aceto, I.~Damg{\aa}rd, L.~A. Goldberg, M.~M. Halld{\'o}rsson,
  A.~Ing{\'o}lfsd{\'o}ttir, and I.~Walukiewicz, editors, {\em Automata,
  Languages and Programming}, volume 5125 of {\em Lecture Notes in Computer
  Science}, pages 698--710, Berlin, 2008. Springer.
\newblock \href {https://doi.org/10.1007/978-3-540-70575-8\_57}
  {\path{doi:10.1007/978-3-540-70575-8\_57}}.

\bibitem{Gluss1962Further}
B.~Gluss.
\newblock Further remarks on line segment curve-fitting using dynamic
  programming.
\newblock {\em Commun. ACM}, 5(8):441--443, 1962.
\newblock \href {https://doi.org/10.1145/368637.368753}
  {\path{doi:10.1145/368637.368753}}.

\bibitem{Golan2003Semirings}
J.~S. Golan.
\newblock {\em Semirings and Affine Equations Over Them}, volume 556 of {\em
  Mathematics and Its Applications}.
\newblock Springer, Dordrecht, 2003.
\newblock \href {https://doi.org/10.1007/978-94-017-0383-3}
  {\path{doi:10.1007/978-94-017-0383-3}}.

\bibitem{Gondran2008Graphs}
M.~Gondran and M.~Minoux.
\newblock {\em Graphs, Dioids and Semirings}, volume~41 of {\em Operations
  Research/ Computer Science Interfaces}.
\newblock Springer, New York, NY, 2008.
\newblock \href {https://doi.org/10.1007/978-0-387-75450-5}
  {\path{doi:10.1007/978-0-387-75450-5}}.

\bibitem{Grigoriev2018Tropical}
D.~Grigoriev.
\newblock Tropical {N}ewton-{P}uiseux polynomials.
\newblock In V.~P. Gerdt, W.~Koepf, W.~M. Seiler, and E.~V. Vorozhtsov,
  editors, {\em Computer Algebra in Scientific Computing}, volume 11077 of {\em
  Lecture Notes in Computer Science}, pages 177--186. Springer, Cham, 2018.
\newblock \href {https://doi.org/10.1007/978-3-319-99639-4\_12}
  {\path{doi:10.1007/978-3-319-99639-4\_12}}.

\bibitem{Hartman1959Onfunctions}
P.~Hartman.
\newblock On functions representable as a difference of convex functions.
\newblock {\em Pacific J. Math.}, 9(3):707--713, 1959.
\newblock \href {https://doi.org/10.2140/pjm.1959.9.707}
  {\path{doi:10.2140/pjm.1959.9.707}}.

\bibitem{Heidergott2006Maxplus}
B.~Heidergott, G.~J. Olsder, and J.~{van der Woude}.
\newblock {\em Max {P}lus at Work}.
\newblock Princeton Series in Applied Mathematics. Princeton Univ. Press,
  Princeton, NJ, 2006.

\bibitem{Imai1986Optimal}
H.~Imai and M.~Iri.
\newblock An optimal algorithm for approximating a piecewise linear function.
\newblock {\em J. Inf. Process.}, 9(3):159--162, 1986.

\bibitem{Itenberg2007Tropical}
I.~Itenberg, G.~Mikhalkin, and E.~Shustin.
\newblock {\em Tropical Algebraic Geometry}, volume~35 of {\em Oberwolfach
  Seminars}.
\newblock Birkh\"{a}user, Basel, 2007.
\newblock \href {https://doi.org/10.1007/978-3-7643-8310-7}
  {\path{doi:10.1007/978-3-7643-8310-7}}.

\bibitem{Kenoufi2025Idempotent}
A.~O. Kenoufi, M.~Gondran, and A.~Gondran.
\newblock {\em Tropical Mathematics and Applications to Theoretical Physics and
  Scientific Computing}.
\newblock De Gruyter, Berlin, 2025.
\newblock \href {https://doi.org/doi:10.1515/9783110769265}
  {\path{doi:doi:10.1515/9783110769265}}.

\bibitem{Kolokoltsov1997Idempotent}
V.~N. Kolokoltsov and V.~P. Maslov.
\newblock {\em Idempotent Analysis and Its Applications}, volume 401 of {\em
  Mathematics and Its Applications}.
\newblock Springer, Dordrecht, 1997.
\newblock \href {https://doi.org/10.1007/978-94-015-8901-7}
  {\path{doi:10.1007/978-94-015-8901-7}}.

\bibitem{Ioannis2025Revisiting}
I.~Kordonis and P.~Maragos.
\newblock Revisiting tropical polynomial division: Theory, algorithms, and
  application to neural networks.
\newblock {\em IEEE Trans. Neural Netw. Learn. Syst.}, 36(9):15978--15992,
  2025.
\newblock \href {https://doi.org/10.1109/TNNLS.2025.3570807}
  {\path{doi:10.1109/TNNLS.2025.3570807}}.

\bibitem{Krivulin2012Solution}
N.~Krivulin.
\newblock A solution of a tropical linear vector equation.
\newblock In S.~Yenuri, editor, {\em Advances in Computer Science}, volume~5 of
  {\em Recent Advances in Computer Engineering Series}, pages 244--249. WSEAS
  Press, Athens, 2012.
\newblock \href {https://arxiv.org/abs/1212.6107} {\path{arXiv:1212.6107}}.

\bibitem{Krivulin2021Algebraic}
N.~Krivulin.
\newblock Algebraic solution of tropical polynomial optimization problems.
\newblock {\em Mathematics}, 9(19):2472, 2021.
\newblock URL: \url{https://www.mdpi.com/2227-7390/9/19/2472}, \href
  {https://arxiv.org/abs/2002.03168} {\path{arXiv:2002.03168}}, \href
  {https://doi.org/10.3390/math9192472} {\path{doi:10.3390/math9192472}}.

\bibitem{Krivulin2023Algebraic}
N.~Krivulin.
\newblock Algebraic solution of tropical best approximation problems.
\newblock {\em Mathematics}, 11(18):3949, 2023.
\newblock \href {https://arxiv.org/abs/2308.07210} {\path{arXiv:2308.07210}},
  \href {https://doi.org/10.3390/math11183949}
  {\path{doi:10.3390/math11183949}}.

\bibitem{Krivulin2024Solution}
N.~Krivulin.
\newblock On solution of tropical discrete best approximation problems.
\newblock {\em Soft Comput.}, 28(20):12097--12112, 2024.
\newblock \href {https://arxiv.org/abs/2403.16337} {\path{arXiv:2403.16337}},
  \href {https://doi.org/10.1007/s00500-024-09940-4}
  {\path{doi:10.1007/s00500-024-09940-4}}.

\bibitem{Krivulin2009Onsolution}
N.~K. Krivulin.
\newblock On solution of a class of linear vector equations in idempotent
  algebra.
\newblock {\em Vestnik of Saint Petersburg University. Applied Mathematics},
  (3):63--76, 2009.
\newblock (in Russian).

\bibitem{Krivulin2023Solution}
N.~K. Krivulin.
\newblock On the solution of a two-sided vector equation in tropical algebra.
\newblock {\em Vestnik St. Petersburg Univ. Math.}, 56(2):172--181, 2023.
\newblock \href {https://doi.org/10.1134/S1063454123020103}
  {\path{doi:10.1134/S1063454123020103}}.

\bibitem{Li1992Morphological}
D.~Li.
\newblock Morphological template decomposition with max-polynomials.
\newblock {\em J. Math. Imaging Vision}, 1(3):215--221, 1992.
\newblock \href {https://doi.org/10.1007/BF00129876}
  {\path{doi:10.1007/BF00129876}}.

\bibitem{Zhiwei2025Achieving}
Z.~Li and C.~Wang.
\newblock Achieving sharp upper bounds on the expressive power of neural
  networks via tropical polynomials.
\newblock {\em IEEE Trans. Neural Netw. Learn. Syst.}, 36(2):2931--2945, 2025.
\newblock \href {https://doi.org/10.1109/TNNLS.2024.3350786}
  {\path{doi:10.1109/TNNLS.2024.3350786}}.

\bibitem{Maclagan2015Introduction}
D.~Maclagan and B.~Sturmfels.
\newblock {\em Introduction to Tropical Geometry}, volume 161 of {\em Graduate
  Studies in Mathematics}.
\newblock AMS, Providence, RI, 2015.
\newblock \href {https://doi.org/10.1090/gsm/161} {\path{doi:10.1090/gsm/161}}.

\bibitem{Maragos2021Tropical}
P.~Maragos, V.~Charisopoulos, and E.~Theodosis.
\newblock Tropical geometry and machine learning.
\newblock {\em Proc. IEEE}, 109(5):728--755, 2021.
\newblock \href {https://doi.org/10.1109/JPROC.2021.3065238}
  {\path{doi:10.1109/JPROC.2021.3065238}}.

\bibitem{Markwig2010Field}
T.~Markwig.
\newblock A field of generalised {P}uiseux series for tropical geometry.
\newblock {\em Rend. Sem. Mat. Univ. Politec. Torino}, 68(1):79--82, 2010.

\bibitem{Mhaskar2000Fundamentals}
H.~N. Mhaskar and D.~V. Pai.
\newblock {\em Fundamentals of Approximation Theory}.
\newblock Narosa Publishing House, New Delhi, 2000.

\bibitem{Osborne1967Best}
M.~R. Osborne and G.~A. Watson.
\newblock On the best linear {C}hebyshev approximation.
\newblock {\em Comput. J.}, 10(2):172--177, 1967.
\newblock \href {https://doi.org/10.1093/comjnl/10.2.172}
  {\path{doi:10.1093/comjnl/10.2.172}}.

\bibitem{Saad2021Zerosum}
O.~Saadi.
\newblock {\em Zero-sum repeated games: Accelerated algorithms and tropical
  best-approximation}.
\newblock PhD thesis, Institut Polytechnique de Paris, 2021.

\bibitem{Sposito1976Minimizing}
V.~A. Sposito.
\newblock Minimizing the maximum absolute deviation.
\newblock {\em ACM SIGMAP Bull.}, (20):51--53, 1976.
\newblock \href {https://doi.org/10.1145/1217073.1217077}
  {\path{doi:10.1145/1217073.1217077}}.

\bibitem{Steffens2006History}
K.-G. Steffens.
\newblock {\em The History of Approximation Theory}.
\newblock Birkh\"{a}user, Boston, MA, 2006.
\newblock \href {https://doi.org/10.1007/0-8176-4475-X}
  {\path{doi:10.1007/0-8176-4475-X}}.

\bibitem{Stiefel1960Note}
E.~Stiefel.
\newblock Note on {J}ordan elimination, linear programming and {T}chebycheff
  approximation.
\newblock {\em Numer. Math.}, 2:1--17, 1960.
\newblock \href {https://doi.org/10.1007/BF01386203}
  {\path{doi:10.1007/BF01386203}}.

\bibitem{Stone1961Approximation}
H.~Stone.
\newblock Approximation of curves by line segments.
\newblock {\em Math. Comp.}, 15:40--47, 1961.
\newblock \href {https://doi.org/10.1090/S0025-5718-1961-0119390-6}
  {\path{doi:10.1090/S0025-5718-1961-0119390-6}}.

\bibitem{Szusz2010Linear}
E.~K. Szusz and A.~R. Willms.
\newblock A linear time algorithm for near minimax continuous piecewise linear
  representations of discrete data.
\newblock {\em SIAM J. Sci. Comput.}, 32(5):2584--2602, 2010.
\newblock \href {https://doi.org/10.1137/090769077}
  {\path{doi:10.1137/090769077}}.

\bibitem{Tharwat2010One}
A.~Tharwat and K.~Zimmermann.
\newblock One class of separable optimization problems: Solution method,
  application.
\newblock {\em Optimization}, 59(5):619--625, 2010.
\newblock \href {https://doi.org/10.1080/02331930801954698}
  {\path{doi:10.1080/02331930801954698}}.

\bibitem{Tomek1974Two}
I.~Tomek.
\newblock Two algorithms for piecewise-linear continuous approximation of
  functions of one variable.
\newblock {\em IEEE Trans. Comput.}, 23(4):445--448, 1974.
\newblock \href {https://doi.org/10.1109/T-C.1974.223961}
  {\path{doi:10.1109/T-C.1974.223961}}.

\bibitem{Tuy2016Convex}
H.~Tuy.
\newblock {\em Convex Analysis and Global Optimization}, volume 110 of {\em
  Springer Optimization and Its Applications}.
\newblock Springer, Cham, 2 edition, 2016.
\newblock \href {https://doi.org/10.1007/978-3-319-31484-6}
  {\path{doi:10.1007/978-3-319-31484-6}}.

\bibitem{Wang2023Tropical}
J.~Wang.
\newblock Tropical algebra with high-order matrix for multiple-noise removal.
\newblock {\em J. Low Freq. Noise Vib. Act. Control}, 42(2):898--910, 2023.
\newblock \href {https://doi.org/10.1177/14613484221143348}
  {\path{doi:10.1177/14613484221143348}}.

\bibitem{Watson1970Algorithm}
G.~A. Watson.
\newblock On an algorithm for nonlinear minimax approximation.
\newblock {\em Commun. ACM}, 13(3):160--162, 1970.
\newblock \href {https://doi.org/10.1145/362052.362056}
  {\path{doi:10.1145/362052.362056}}.

\bibitem{Zhang2018Tropical}
L.~Zhang, G.~Naitzat, and L.-H. Lim.
\newblock Tropical geometry of deep neural networks.
\newblock In J.~Dy and A.~Krause, editors, {\em Proceedings of the 35th
  International Conference on Machine Learning}, volume~80 of {\em Proceedings
  of Machine Learning Research}, pages 5824--5832. PMLR, Cambridge, MA, 2018.

\bibitem{Zimmermann1984Maxseparable}
K.~Zimmermann.
\newblock On max-separable optimization problems.
\newblock In R.~E. Burkard, R.~A. Cuninghame-Green, and U.~Zimmermann, editors,
  {\em Algebraic and Combinatorial Methods in Operations Research}, volume~95
  of {\em North-Holland Mathematics Studies}, pages 357--362. North-Holland,
  Amsterdam, 1984.
\newblock \href {https://doi.org/10.1016/S0304-0208(08)72967-0}
  {\path{doi:10.1016/S0304-0208(08)72967-0}}.

\bibitem{Zimmermann2003Disjunctive}
K.~Zimmermann.
\newblock Disjunctive optimization, max-separable problems and extremal
  algebras.
\newblock {\em Theoret. Comput. Sci.}, 293(1):45--54, 2003.
\newblock \href {https://doi.org/10.1016/S0304-3975(02)00231-1}
  {\path{doi:10.1016/S0304-3975(02)00231-1}}.

\end{thebibliography}

\end{document}